\documentclass[12pt]{amsart}

\usepackage{epsfig}
\usepackage{mathtools}
\usepackage[subnum]{cases}

\usepackage{amsmath,amssymb}
\usepackage{color}
\newtheorem{theorem}{Theorem}[section]
\newtheorem{definition}[theorem]{Definition}
\newtheorem{lemma}[theorem]{Lemma}
\newtheorem{corollary}[theorem]{Corollary}
\newtheorem{remark}[theorem]{Remark}
\newtheorem{proposition}[theorem]{Proposition}

\newtheorem{problem}[theorem]{Problem}

\newcommand{\inte}{{\mathrm{int}}\,}
\newcommand{\reg}{{\mathrm{reg}}\,}
\newcommand{\relint}{{\mathrm{relint}}\,}
\newcommand{\cl}{{\mathrm{cl}}\,}
\newcommand{\conv}{{\mathrm{conv}}\,}

\newcommand{\ee}{\varepsilon}

\newcommand{\K}{\mathcal{K}}

\def\cK{\mathcal{K}}
\def\cH{\mathcal{K}}
\def\sphere{S^{n-1}}
\def\N{\mathbb{N}}
\def\Rn{{\mathbb R^n}}
\def\R{\mathbb{R}}
\def\cH{\mathcal{H}}

\def\deV{\widetilde{C}_{H, \psi}}

\def\leV{\widetilde{C}_{G}}
\def\deV{\widetilde{C}_{G, \psi}}
\def\deVo
{\widetilde{C}_{G, \psi_0}}

\def\ball{B^n}

\def\wp{G}
\def\dveV{\widetilde{V}_G}

\def\bt{\begin{theorem}}
	\def\et{\end{theorem}}
\def\bl{\begin{lemma}}
	\def\el{\end{lemma}}
\def\br{\begin{remark}}
	\def\er{\end{remark}}
\def\bc{\begin{corollary}}
	\def\ec{\end{corollary}}
\def\bd{\begin{definition}}
	\def\ed{\end{definition}}
\def\bp{\begin{proposition}}
	\def\ep{\end{proposition}}

\begin{document}

\title{GENERAL VOLUMES IN THE ORLICZ-BRUNN-MINKOWSKI THEORY AND A RELATED MINKOWSKI PROBLEM II}
\author[Richard J. Gardner, Daniel Hug, Sudan Xing, and Deping Ye]
	{Richard J. Gardner, Daniel Hug, Sudan Xing, and Deping Ye}
\address{Department of Mathematics, Western Washington University,
		Bellingham, WA 98225-9063, USA} \email{richard.gardner@wwu.edu}
\address{Karlsruhe Institute of Technology, Department of Mathematics,
		D-76128 Karlsruhe, Germany}	
\email{daniel.hug@kit.edu}
	\address{Department of Mathematics and Statistics, Memorial University of Newfoundland,\newline St.~John's, Newfoundland, Canada A1C 5S7} \email{sudanxing@gmail.com}
	\address{Department of Mathematics and Statistics, Memorial University of Newfoundland,\newline St.~John's, Newfoundland, Canada A1C 5S7} \email{deping.ye@mun.ca}
\thanks{First author supported in part by U.S.~National Science Foundation Grant DMS-1402929.  Second author supported in part by German Research Foundation (DFG) grants HU 1874/4-2 and FOR 1548.  Fourth author supported in part by an NSERC grant.}	

\begin{abstract}
The general dual volume $\dveV(K)$ and the general dual Orlicz curvature measure $\deV(K, \cdot)$ were recently introduced for functions $G: (0, \infty)\times \sphere\rightarrow (0, \infty)$ and convex bodies $K$ in $\R^n$ containing the origin in their interiors.  We  extend $\dveV(K)$ and $\deV(K, \cdot)$ to more general functions  $G: [0, \infty)\times \sphere\rightarrow [0, \infty)$ and to compact convex sets $K$ containing the origin (but not necessarily in their interiors). Some basic properties of the general dual volume and of the dual Orlicz curvature measure, such as the continuous dependence on the underlying set, are provided. These are required to study a Minkowski-type problem for the dual Orlicz curvature measure. We mainly focus on the case when $G$ and $\psi$ are both increasing, thus complementing our previous work.

The Minkowski problem asks to characterize Borel measures $\mu$ on $\sphere$ for which there is a convex body $K$ in $\R^n$ containing the origin such that $\mu$ equals $\deV(K, \cdot)$, up to a constant.  A major step in the analysis concerns discrete measures $\mu$, for which we prove the existence of convex polytopes containing the origin in their interiors solving the Minkowski problem. For general (not necessarily discrete) measures $\mu$, we use an approximation argument. This approach is also applied to the case where $G$ is decreasing and $\psi$ is increasing, and hence augments our previous work. When the measures $\mu$ are even, solutions that are origin-symmetric convex bodies are also provided under some mild conditions on $G$ and $\psi$. Our results generalize several previous works and provide more precise information about the solutions of the Minkowski problem when $\mu$ is discrete or even.
\end{abstract}

\maketitle
\pagestyle{myheadings}
\markboth{RICHARD J. GARDNER, DANIEL HUG, SUDAN XING, AND DEPING YE}{THE ORLICZ-BRUNN-MINKOWSKI THEORY AND A MINKOWSKI PROBLEM}

\section{Introduction}\label{intro}
The study of Minkowski problems, initiated by Minkowski \cite{min1897, min1903} over a century ago, took on a new life when Lutwak  \cite{Lu93} introduced the $L_p$ surface area measures for $p\ge 1$ and the corresponding $L_p$ Minkowski problem, the case $p=1$ of which is the classical Minkowski problem.  For $p\in \R$, the $L_p$ Minkowski problem asks for necessary and/or sufficient conditions for a measure $\mu$ on the unit sphere $\sphere$ in $\R^n$ to be the $L_p$ surface area measure of a convex body (i.e., a convex compact set in $\Rn$ with nonempty interior).  The $L_p$ surface area measures can be obtained via a first-order variation of volume with respect to $L_p$ addition of convex bodies, making the method of Lagrange multipliers a powerful tool in solving the $L_p$ Minkowski problem; see e.g.~\cite{HugLYZ, Lu93, LYZ04}. The case $p=0$ is of particular significance because the $L_0$ surface area measure, the so-called cone volume measure, is an affine invariant. The $L_0$ or logarithmic Minkowski problem is challenging and was only solved recently for even measures by B\"{o}r\"{o}czky, Lutwak, Yang, and Zhang \cite{BLYZ2013}.  More recent contributions to the logarithmic Minkowski problem are \cite{BHZ2016, zhug2014} and further references and background on the $L_p$ Minkowski problem may be found in \cite{LYZActa, Sch}.

Recent seminal work of Huang, Lutwak, Yang, and Zhang \cite{LYZActa} brought new ingredients, the $q$th dual curvature measures, to the family of Minkowski problems.   These measures are obtained via a first-order variation of the $q$th dual volume with respect to $L_0$ addition of convex bodies (see \cite[Theorem 4.5]{LYZActa}), the case $q=n$ being the $L_0$ surface area. The authors of \cite{LYZActa} posed a corresponding Minkowski problem---the dual Minkowski problem---of finding necessary and/or sufficient conditions for a measure $\mu$ on $\sphere$ to be the $q$th dual curvature measure of some convex body, and they provided a partial solution when $\mu$ is even.  Clearly, the logarithmic Minkowski problem is a special case of the dual Minkowski problem.  Naturally, the dual Minkowski problem has become important for the dual Brunn-Minkowski theory introduced by Lutwak \cite{Lu75, Lu88}.  Since \cite{LYZActa}, progress includes a complete solution for $q<0$ by Zhao \cite{zhao}, solutions for even $\mu$ in \cite{BHP, BLYZZ2017, HP-adv,  zhao-jdg}, and solutions via curvature flows and partial differential equations in \cite{ChenLi, JiangWu, LiShengWang}.

An important extension of the dual Minkowski problem was carried out by Lutwak, Yang, and Zhang \cite{LYZ-Lp}, who introduced $L_p$ dual curvature measures and posed a corresponding $L_p$ dual Minkowski problem.  In \cite{LYZ-Lp}, the $L_0$ addition in \cite{LYZActa} is replaced by $L_p$ addition, while the $q$th dual volume remains unchanged.  The first contribution to the $L_p$ dual Minkowski problem, by Huang and Zhao \cite{HuangZhao}, proves the existence of solutions for $p, q\in\R$ when $p>0$ and $q<0$, and for even $\mu$ when $pq>0$, $p\neq q$. Their results were augmented by Chen, Huang, and Zhao \cite{Chen-H-Z}, who used curvature flows to show the smoothness of solutions for even $\mu$ and $pq \geq 0$. B\"{o}r\"{o}czky and Fodor \cite{BorFor} provide a beautiful solution to the $L_p$ dual Minkowski problem for general $\mu$ when $p>1$ and $q>0$.

The first Orlicz version of the Minkowski problem appeared in \cite{HLYZ2010-adv}, at the inception of the Orlicz-Brunn-Minkowski theory in 2010.  Later, the dual Minkowski problem was extended to the Orlicz setting and partially solved in \cite{XY2017-1, ZSY2017}; here, the $q$th dual volume in \cite{LYZActa} is replaced by certain dual Orlicz quermassintegrals, while $L_0$ addition is retained.  A common generalization of the problems in \cite{LYZ-Lp, XY2017-1, ZSY2017} was proposed in \cite{GHWXY}, in which the $q$th dual volume is replaced by a very general notion of dual volume denoted by $\dveV(K)$ and, simultaneously, $L_0$ addition is replaced by an extension of $L_p$ addition called Orlicz addition.   By combining the general dual volume $\dveV(K)$ with Orlicz addition, a general dual Orlicz curvature measure denoted by $\deV(K, \cdot)$ is defined in \cite{GHWXY}, where $K$ is a convex body with the origin in its interior, and $G: (0, \infty)\times \sphere\rightarrow (0, \infty)$ and $\psi: (0, \infty)\rightarrow (0, \infty)$ are continuous.  The two-variable function $G$ allows $\dveV(K)$ to include not only the $q$th dual volume and the dual Orlicz quermassintegrals of \cite{XY2017-1, ZSY2017}, but several other related notions as well.  In \cite{GHWXY}, the following Minkowski problem (see Problem \ref{Minkowski-c-06-25} below) was stated:

{\em For which nonzero finite Borel measures $\mu$ on $\sphere$ and continuous functions $\wp:(0, \infty)\times \sphere\rightarrow (0, \infty)$ and $\psi:(0,\infty)\to(0,\infty)$ do there exist $\tau\in \R$ and a convex body $K$ (with the origin in its interior) such that $\mu=\tau\,\deV(K,\cdot)$?}

The problem, which requires solving a certain Monge-Amp\`{e}re equation (see (\ref{MPnew}) below), contains all previously known Minkowski problems as special cases. A solution was presented in \cite[Theorem~6.4]{GHWXY} for general measures $\mu$, assuming that $G_t=\partial G(t, u)/\partial t<0$, $G$ satisfies some growth conditions, and  $\psi$ satisfies (\ref{June-22}) below.  Our aim here is to complement our results in \cite{GHWXY} by dealing with the case when $G_t>0$. This requires extending $\dveV(K)$ and $\deV(K, \cdot)$ to more general functions $G: [0, \infty)\times \sphere\rightarrow [0, \infty)$ and to compact convex sets $K$ containing the origin, but not necessarily in their interiors; see Sections~\ref{section-3} and~\ref{section-5}, especially Definitions~\ref{july251} and~\ref{general dual Orlicz curvature measure-7-01}. It is also necessary to show that $\dveV(K)$ and $\deV(K, \cdot)$ are continuous in $K$ (see Lemma~\ref{continuity for vG} and Proposition~\ref{property-extended-7-1}(iii)), a task necessitating a more delicate treatment of the various maps and cones related to a compact convex set (see Sections \ref{Maps} and \ref{section-5}) than when the origin is contained in the interior.

Unlike the proof of \cite[Theorem~6.4]{GHWXY}, we approach the Minkowski problem stated above when $G_t>0$ by first dealing with the case when $\mu$ is discrete. This is achieved in Theorem \ref{whole-solution-6-27}, where we establish, under certain growth conditions on $\psi$, the existence of a convex polytope $P$ with the origin in its interior, such that $\mu$ equals $\deV(P, \cdot)$ (up to a normalization constant).  If $G(t, u)=t^n/n$, then $\dveV(K)$ is the volume of $K$, so our Theorem  \ref{whole-solution-6-27} recovers the solutions to the Orlicz-Minkowski problem for discrete measures by Huang and He \cite{huang} and Li \cite{liaijun2014}. When $\psi(t)=t^p$ for $p>1$ and $G(t, u)=t^q\phi(u)$ for $q>0$ and $\phi\in C^+(\sphere)$, Theorem  \ref{whole-solution-6-27} recovers the solution to the $L_p$ dual Minkowski problem for discrete measures by B\"{o}r\"{o}czky and Fodor \cite[Theorem~1.1]{BorFor}.  The techniques in these works are similar and based on those in \cite{HugLYZ}, but some of our arguments differ from and are rather more complicated than those in \cite{BorFor, huang, liaijun2014}.  In particular, the general volume $\dveV(\cdot)$ prohibits the use of Minkowski's inequality as in \cite{huang, liaijun2014}, and in general the two-variable function $G$, and the lack of homogeneity of $G$ and $\psi$, require a somewhat more delicate analysis than the special case considered in  \cite{BorFor}.   On the other hand, we are able to avoid some constructions in \cite{BorFor} by making use of the absolute continuity of $\deV(K,\cdot)$ with respect to the surface area measure of $K$ proved in Proposition~\ref{property-extended-7-1}(ii).

With Theorem~\ref{whole-solution-6-27} in hand, our Minkowski problem for general measures $\mu$ can be solved by approximation.  This is accomplished in Theorem~\ref{whole-general-solution-7-5}, where it is shown that under certain conditions on $G$ and $\psi$, including $G_t>0$, a finite Borel measure $\mu$ on $\sphere$ is not concentrated on any closed hemisphere if and only if there exists a convex body $K$ containing the origin such that
\begin{equation*}
(\psi\circ h_K)\mu= \bigg(\int_{\sphere} \psi(h_{K}(u))\,d\mu(u)\bigg)  \ \frac{\leV(K, \cdot)}{\leV(K, \sphere)},
\end{equation*}
where $h_K$ denotes the support function of $K$ (see Section~\ref{prel} for notation and most definitions).  Again, this result recovers (in a slightly different form) and strengthens the solutions to the Orlicz-Minkowski problem in \cite[Theorem~1.2]{huang} and the  $L_p$ dual Minkowski problem in  \cite[Theorem~1.2]{BorFor}.  In Theorem~\ref{whole-general-solution-7-9-decreaing}, we use the same approximation technique to prove a variant of \cite[Theorem~6.4]{GHWXY} in the case when $G_t<0$.  When $\psi(t)=t^p$, $p>0$, and $G(t, u)=t^q$, $q<0$, Theorem~\ref{whole-general-solution-7-9-decreaing} implies \cite[Theorem~3.5]{HuangZhao}.  We end Section~\ref{section-6} with Theorem~\ref{uniquemp}, a uniqueness result related to Theorem~\ref{whole-general-solution-7-9-decreaing} under some additional assumptions on the underlying convex bodies.  As far as we know, this is the first uniqueness result for Orlicz-Minkowski problems that applies when $G(t, u)=t^n/n$ and $\dveV(K)$ is the volume of $K$.  A special case of Theorem~\ref{uniquemp} contributes to \cite[Problem~8.2]{LYZ-Lp} by providing a counterpart to \cite[Theorem~8.3]{LYZ-Lp} for sufficiently smooth convex bodies and generalizing the uniqueness assertion in \cite[Theorem~4.1]{HuangZhao}.  The uniqueness problem for general dual Orlicz curvature measures remains open.

In Section~\ref{section-7}, we focus on the case when the measure $\mu$ is even, in which case one expects the solution to be an origin-symmetric convex body. Each such body generates a norm on $\Rn$, and every norm on $\R^n$ arises from an origin-symmetric convex body.  This lends special significance to Minkowski problems for even measures, particularly in applications to analysis; for example, in proving the $L_p$ affine Sobolev inequality \cite{LYZ02} and the affine Moser-Trudinger and Morrey-Sobolev inequalities \cite{CLYZ-1}.  Corresponding to Theorems~\ref{whole-general-solution-7-5} and ~\ref{whole-general-solution-7-9-decreaing}, we prove Theorems~\ref{whole-general-solution-7-55-even} and~\ref{whole-general-solution-7-9-decreaing-even} for even $\mu$, where it is natural to impose weaker conditions on $\psi$ but an extra assumption on $G$ (i.e., that $G_t(t, \cdot)$ is even in $t$). In our final result, Theorem~\ref{solution-general-dual-Orlicz-main theorem-11-27-even}, we solve our Minkowski problem under the assumption that $\mu$ is an even measure vanishing on any great subsphere, when $G_t<0$ and $\psi$ is decreasing.  Again, if $G(t, u)=t^n/n$, $\dveV(K)$ is the volume of $K$ and we recover the solution to the Orlicz-Minkowski problem for even measures by Haberl, Lutwak, Yang, and Zhang \cite{HLYZ2010-adv}.  Moreover, when $\psi(t)=t^p$ and $G(t, u)=t^q$,  Theorems~\ref{whole-general-solution-7-55-even} and and~\ref{solution-general-dual-Orlicz-main theorem-11-27-even} yield the results of Huang and Zhao \cite[Theorem~3.9]{HuangZhao} for $p, q>0$ and $p\neq q$ and \cite[Theorem~3.11]{HuangZhao} for $p, q<0$ and $p\neq q$, respectively. The method we employ both avoids the use of John ellipsoids in the proof of \cite[Theorem~3.9]{HuangZhao} and provides detailed information, not given in \cite{HuangZhao}, on the polytopal solutions to our Minkowski problem when $\mu$ is an even discrete measure.

\section{Preliminaries and Background}\label{prel}
We work in $\R^n$ with scalar product $\langle \cdot, \cdot \rangle$ and  Euclidean norm $|\cdot|$.  The origin and canonical orthonormal basis are denoted by $o$ and $\{e_1, \ldots, e_n\}$, respectively.  Let $\ball=\{x\in \Rn: |x|\leq 1\}$ and $\sphere=\{x\in \Rn: |x|=1\}$ be the unit ball and unit sphere in $\Rn$.  The characteristic function of a set $E$ is signified by $1_E$.

If $X$ is a set,  we denote by $\conv X$, $\cl X$, $\inte X$, $\relint X$, $\partial X$, and $\dim X$ the {\it convex hull}, {\it closure}, {\it interior},
{\it relative interior} (that is, the interior with respect to the affine hull), {\it boundary}, and {\it dimension} (that is, the dimension of the affine hull) of $X$, respectively.  We say that $X$ is {\em origin symmetric} if it is invariant under reflection in $o$.  If $x\in \R^n\setminus\{o\}$, then $x^{\perp}$ is the $(n-1)$-dimensional subspace orthogonal to $x$.

We write ${\mathcal{H}}^k$ for $k$-dimensional Hausdorff measure in $\R^n$, where $k\in \{1,\dots,n\}$.  For compact sets $E$, we also write $V_n(E)={\mathcal{H}}^n(E)$ for the {\em volume} of $E$.  The volume of the unit ball is $\kappa_n=V_n(B^n)$ and then ${\mathcal{H}}^{n-1}(\sphere)=n\kappa_n$.  The notation $dx$ means $d{\mathcal{H}}^k(x)$ for the appropriate $k\in\{1,\dots,n\}$, unless stated otherwise. In particular, integration on $\sphere$ is usually denoted by $du=d{\mathcal{H}}^{n-1}(u)$.

The class of nonempty compact convex sets in $\R^n$ is written $\cK^n$. We will often work with $\cK_o^n$, the members of $\cK^n$ containing $o$, or $\cK_{(o)}^n$, the {\em convex bodies} (i.e., compact convex subsets in $\Rn$ with nonempty interiors) {\em containing $o$ in their interiors}. A {\em convex polytope} is the convex hull of a finite subset of $\Rn$.  For the following information about convex sets, we refer the reader to \cite{Gruber2007, Sch}. The standard metric on $\cK^n$ is the {\em Hausdorff metric} $\delta(\cdot, \cdot)$, which can be defined by
$$\delta(K, L)=\|h_K-h_ {L}\|_{\infty} =\sup_{u\in \sphere} |h_K(u)-h_{L}(u)|$$
for $K, L\in \cK^n$, where $h_K: \sphere\rightarrow \R$ is the {\em support function} of $K\in \cK^n$, given by $h_K(u)=\max\{\langle u, x\rangle: x\in K\}$ for $u\in \sphere$. Note that
$$H(K,u)=\big\{x\in\R^n: \langle x, u\rangle=h_K(u)\big\}$$
is the supporting hyperplane of $K$ in the direction $u\in \sphere$. The corresponding support set of $K$ in direction $u$ is $F(K,u)=K\cap H(K,u)$.
We say that the sequence $K_1, K_2,\ldots$ of sets in $\cK^n$ {\em converges} to $K\in \cK^n$ if $\lim_{i\rightarrow \infty} \delta(K_i, K)=0$.  The {\em Blaschke selection theorem} states that every bounded sequence in $\cK^n$ has a subsequence that converges to a set in $\cK^n$.

The {\em surface area measure} $S(K,\cdot)$ of a convex body $K$ in $\R^n$ is defined for Borel sets $E\subset \sphere$ and satisfies
$$
S(K,E)=\mathcal{H}^{n-1}(\nu_{K}^{-1}(E)),
$$
where $\nu_{K}^{-1}(E)=\{x\in \partial K: x\in H(K,u) \text{ for some }u\in E\}$
 (see Section~\ref{Maps}).

Let $\mu$ be a nonzero finite Borel measure on $\sphere$.  We say that $\mu$ is {\em not concentrated on any closed hemisphere}  if \begin{equation}\label{condition for Minkowski problem}
\int _{\sphere} \langle u, v\rangle_+ \,d\mu(u)>0\ \ \ \mathrm{for~} v\in \sphere,
\end{equation}
where  $a_+=\max\{a, 0\}$ for $a\in \R$.

As usual, $C(E)$ denotes the class of continuous functions on a topological space $E$, and we shall write $C^+(E)$ for the (strictly) positive functions in $C(E)$. Let $\Omega\subset  \sphere$ be a closed set not contained in any closed hemisphere of $\sphere$.  For each $f\in C^+(\Omega)$, one can define a convex body $[f]$, the {\em Aleksandrov body} (or {\em Wulff shape}) associated to it, by setting
\begin{equation}\label{july287}
[f]= \bigcap_{u\in\Omega}\big\{x\in\R^n:\langle x, u\rangle\leq f(u)\big\}.
\end{equation}
In particular, when $\Omega=\sphere$ and $f=h_K$ for $K\in \cK^n$, one has
$$K=[h_K]=\bigcap_{u\in\sphere}
\big\{x\in\R^n:\langle x, u\rangle\leq h_K(u)\big\}.$$

For $K\in \K^n_o$ and  $x\in \Rn\setminus\{o\}$, the {\em radial function} $\rho_K:\Rn\setminus \{o\}\rightarrow [0,\infty)$ of $K$ is defined by
$$\rho_K(x)=\max\{\lambda\ge 0: \lambda x\in K\}.$$
The function $\rho_K$ is homogeneous of degree $-1$, that is, $\rho_K(rx)=r^{-1}\rho_K(x)$ for $x\in \R^n\setminus \{o\}$.

If $K\in\mathcal{K}_{(o)}^n$, the {\em polar body} $K^*$ of $K$ is defined by
$$K^*=\{x\in\mathbb{R}^n: \langle x, y\rangle\leq1 \ \text{for}\  y\in K \}.$$
Then $(K^*)^*=K$ and (see \cite[(1.52), p.~57]{Sch})
\begin{eqnarray} \label{bi-polar--1}
\rho_K(x) h_{K^*}(x)=h_K (x) \rho_{K^*}(x)=1\ \ \ \mathrm{for~} x\in \R^n\setminus\{o\}.
\end{eqnarray}

For $a\in \R\cup\{-\infty\}$, we require the following class of functions $\varphi:(0,\infty)\to (a,\infty)$:
\begin{equation}\label{aug221}
{\mathcal{J}}_a=\{\mbox{$\varphi$ is continuous and strictly monotonic,   $\inf_{t>0}\varphi(t)=a$, and  $\sup_{t>0}\varphi(t)=\infty$}\}.
\end{equation}

Let $h_0\in C^+(\sphere)$, let $g\in C(\sphere)$, and let $\varphi_0\in {\mathcal{J}}_a$ for some $a\in \R\cup\{-\infty\}$. Then we have $\varphi_0^{-1}: (a,\infty)\to (0,\infty)$, and since $\sphere$ is compact, $0<c\le h_0\le C$ for some $0<c\le C$.  It is then easy to check that for $\ee\in \R$ close to $0$, one can define $h_{\ee}=h_{\ee}(h_0,g,\varphi_0)\in C^+(\sphere)$ by
\begin{eqnarray}\label{genplus}
h_{\ee}(u)= {\varphi_0}^{-1}\left(\varphi_0(h_0(u))+\ee g(u)\right).
\end{eqnarray}
In particular, we can apply (\ref{genplus}) when $h_0=h_K$ for some $K\in \cK^n_{(o)}$.

\subsection{Maps and cones related to a compact convex set}\label{Maps}
For $K\in\cK^n$, the {\em normal cone} of $K$ at $z\in K$ is defined by
$$N(K, z)=\{y\in\R^n:\,\langle y,x-z\rangle\leq 0\mbox{ \ for $x\in K$}\}.$$
This is a closed convex cone, and $N(K, z)=\{o\}$ if $z\in \inte K$.  In particular, if $o\in \partial K$, then
\begin{equation}\label{july201}
N(K,o)=\{y\in\R^n:\,\langle y,x\rangle\leq 0\mbox{ \ for $x\in K$}\}.
\end{equation}

Let $K\in\K^n_o$. Then the {\em dual cone} $N(K, o)^*$ of $N(K, o)$ is given by
\begin{equation}\label{july135}
N(K, o)^*=\{x\in\R^n:\langle x,y\rangle\le 0 \mbox{ \ for $y\in N(K,o)$}\}=\cl\{\lambda x: x\in K \ \mathrm{and}\ \lambda\geq 0\};
\end{equation}
the set on the right side is called the {\em support cone} of $K$ at $o$ (see \cite[p.~91]{Sch}).
If $o\in \inte K$, then  $N(K,o)^*=\R^n$, and if $o\in\partial K$, then
\begin{equation}\label{dual-cone-6-26}
\cH^{n-1}(\sphere\cap \partial N(K, o)^*)=0.
\end{equation}

For $K\in\K^n_o$, let $r_K:\sphere\to\partial K$ denote the {\em radial map} of $K$, defined by $r_K(u)=\rho_K(u)u$.
If $o\in\partial K$, then $r_K$ need not be injective,  since
$\rho_K(u)u=o$ for $u\in \sphere \setminus N(K, o)^*$. The radial map also need not be continuous, but it is upper semicontinuous and hence Borel measurable. However, the restriction of the radial map to $\sphere\cap \relint N(K, o)^*$ is injective
and locally bi-lipschitz. Moreover,
\begin{equation}\label{july221}
\rho_K(u) \, \begin{cases} =0& \text{ if } u\in \sphere \setminus N(K, o)^*,\\
>0& \text{ if } u\in\sphere\cap \relint N(K,o)^*.
\end{cases}
\end{equation}
For $u\in \sphere\cap \partial N(K,o)^*$ we only have $\rho_K(u)\ge 0$, but if $K$ is a convex polytope, then
\begin{equation}\label{aug231}
\rho_K(u)>0~~\quad{\text{if and only if}}~~\quad u\in \sphere \cap N(K, o)^*.
\end{equation}

We recall some terminology and facts from \cite[Section~2.2]{LYZActa} and \cite[Section~2.2]{Sch},
presented in a slightly different form (see also \cite{BorFor}).
The {\em radial projection} $\tilde{\pi}: \Rn\setminus \{o\}\to \sphere$ is defined by $\tilde{\pi}(x)=x/|x|$.
For $K\in\K^n$ and $E\subset \partial K$, the {\em spherical image} of $E$ is defined by
$$
\pmb{\nu}_K(E)=\sphere\cap \bigcup_{x\in E} N(K,x).
$$
Recall that for a Borel set $E\subset\partial K$, the spherical image $\pmb{\nu}_K(E)\subset \sphere$ is $\cH^{n-1}$-measurable (see \cite[Lemma 2.2.13]{Sch}).

For the remainder of this section, assume that $K\in\cK_o^n$. Following \cite[p.~340]{LYZActa}, the {\em radial spherical image} of a set $E \subset \sphere$ is given by
$$\pmb{\alpha}_K(E)=\pmb{\nu}_K(r_K(E))\subset\sphere.$$
If $E\subset \sphere\cap \relint N(K,o)^*$ is a Borel set, then so is $r_K(E)$, so $\pmb{\alpha}_K(E)$ is $\cH^{n-1}$-measurable. If $\emptyset\neq E\subset\sphere\setminus N(K,o)^*$, then $r_K(E)=\{o\}$ and again $\pmb{\alpha}_K(E)\subset\sphere\cap N(K,o)$ is $\cH^{n-1}$-measurable. The situation for a Borel set $E\subset \sphere$ contained in the relative boundary of $ N(K,o)^*$ seems to be unclear but will not be needed.

If $\dim K=n$ and $\sigma_K\subset\partial K$ and $\omega_K\subset \sphere\cap  N(K, o)^*$ are the sets of $x\in\partial K$ and $u\in \sphere\cap  N(K, o)^*$ such that $\pmb{\nu}_K(\{x\})$ and $\pmb{\alpha}_K(\{u\})$, respectively, have two or more elements, then
\begin{equation}\label{zerosets}
\cH^{n-1}(\sigma_K)=\cH^{n-1}(\omega_K)=0.
\end{equation}
The sets $\sigma_K$ and $\omega_K$ are $\cH^{n-1}$-measurable, since $\sigma_K$ is a countable union of closed sets and in view of \eqref{dual-cone-6-26}.  The set of {\em regular boundary points} of $K$ is $\reg K=\partial K\setminus\sigma_K$.  We write $\nu_K(x)$ and $\alpha_K(u)$ instead of $\pmb{\nu}_K(\{x\})$ and $\pmb{\alpha}_K(\{u\})$ if $x\in \reg K$ and $u\in \sphere\cap  N(K, o)^* \setminus \omega_K$, respectively. The map $\nu_K$ is continuous on $\reg K$ and $\alpha_K$ is continuous on $\sphere\cap \relint N(K,o)^*\setminus \omega_K$. Hence $\nu_K$ is measurable with respect to $\cH^{n-1}$ on $\partial K$ (if extended arbitrarily on $\partial K\setminus\reg K$) and $\alpha_K$ is measurable with respect to $\cH^{n-1}$ on $\sphere\cap  N(K,o)^*$ (if extended arbitrarily).

The {\em reverse spherical image} of $E\subset \sphere$ is given by
$$\pmb{\nu}_K^{-1}(E)=\{x\in \partial K :N(K,x)\cap E\neq \emptyset \}.$$
If $E$ is a Borel set, then $\pmb{\nu}_K^{-1}(E)$ is $\cH^{n-1}$-measurable. Let $\dim K=n$. Then the set of $u\in\sphere$  such that $\pmb{\nu}_K^{-1}(\{u\})$ has at least two elements is called the set of {\em singular normal vectors} and is a Borel set of $\cH^{n-1}$-measure zero.
Note that $\pmb{\nu}_K$ and $\pmb{\nu}^{-1}_K$ do not necessarily map disjoint sets to disjoint sets. However, the intersections are sets of singular normal vectors and sets of singular boundary points, respectively, and hence have $\cH^{n-1}$-measure zero.  To simplify notation, we will simply write $\nu_K^{-1}$ instead of $\pmb{\nu}_K^{-1}$, since if $u\in\sphere$ is not a singular normal vector and $\nu^{-1}_K(u)$ is the unique point in $\partial K$ with outer unit normal vector $u$, then $\pmb{\nu}^{-1}_K(\{u\})=\{\nu^{-1}_K(u)\}$.

For later use, we also define
\begin{equation}\label{july252}
\Xi_K=\pmb{\nu}_K^{-1}\left(\sphere\cap N(K,o)\right)=K\cap \partial N(K,o)^*.
\end{equation}
Clearly, $\Xi_K=\emptyset$ if $o\in\inte K$. Moreover, if $\dim K\le n-1$, then $\Xi_K=K$.

Following \cite[p.~340]{LYZActa} again, we define the {\em reverse radial spherical image} of a set $E\subset\sphere$ by
\begin{equation}\label{july151}
\pmb{\alpha}^*_K(E)= r_K^{-1}\left(\pmb{\nu}_K^{-1}\left(E\right)\right)\subset \sphere.
\end{equation}
If $E\subset\sphere$ is a Borel set, then $\pmb{\alpha}^*_K(E)$ is $\cH^{n-1}$-measurable. This is shown in \cite[Lemma~2.1]{LYZActa} when $o\in\inte K$. To see that it is true in general, first observe that $\pmb{\nu}_K^{-1}(E)\subset\partial K$ is $\cH^{n-1}$-measurable. If $A=\pmb{\nu}_K^{-1}(E)\cap \relint N(K,o)^*$, then since $r_K$ is locally bi-lipschitz on $\relint N(K,o)^*$, it follows that $r_K^{-1}(A)$ is also $\cH^{n-1}$-measurable. Let $B$ denote the intersection of $\pmb{\nu}_K^{-1}(E)$ with the relative boundary of $N(K,o)^*$. Then $r_K^{-1}(B)\subset
\sphere\cap\partial N(K,o)^*$, which is $\cH^{n-1}$-measurable due to \eqref{dual-cone-6-26}. Therefore $\pmb{\alpha}^*_K(E)=r_K^{-1}(A\cup B)$ is $\cH^{n-1}$-measurable.

For  $u\in \sphere\cap \inte N(K,o)^*\setminus\omega_K$, and hence for  $\cH^{n-1}$-almost all  $u\in \sphere\cap \inte N(K,o)^*$, we have
\begin{equation} \label{relation-11-27}
u\in \pmb{\alpha}_K^*(E)\quad \text{if and only if}\quad  \alpha_K(u)\in E.
\end{equation}

Finally, we remark that
$$
\pmb{\alpha}_K^*\left(E\right)\cap \relint N(K,o)^*\subset  \pmb{\alpha}_K^*\left(E\setminus N(K,o) \right)\subset   \pmb{\alpha}_K^*\left(E\right)\cap   N(K,o)^*,
$$
Examples show that both inclusions can be strict, but in view of \eqref{dual-cone-6-26}, we have
\begin{equation}\label{aug232}
\pmb{\alpha}_K^*\left(E\right)\cap \relint N(K,o)^*= \pmb{\alpha}_K^*\left(E\setminus N(K,o) \right)=   \pmb{\alpha}_K^*\left(E\right)\cap   N(K,o)^*
\end{equation}
up to sets of $\cH^{n-1}$-measure zero.
	
\section{The general dual volume on compact convex sets}\label{section-3}
We begin with the following definition.

\bd\label{july251}
Let $\wp:[0, \infty)\times \sphere\to [0,\infty)$ be such that $u\mapsto\wp(\rho_K(u),u)$ is integrable on $\sphere$ for $K\in \mathcal{K}^n_{o}$. Define the  {\em general dual volume} $\dveV(K)$ of $K\in \mathcal{K}^n_{o}$ by
\begin{equation}\label{def-H-volume-06-23}
\dveV(K)=\int_{\sphere} \wp(\rho_K(u), u)\,du.
\end{equation}
\ed

If $K\in \cK_o^n$ has empty interior, then $\rho_K=0$ outside a great subsphere of $\sphere$.  Since ${\mathcal{H}}^{n-1}$ vanishes on such great subspheres, we then have
\begin{equation}\label{july132}
\dveV(K)=\int_{\sphere}G(0, u)\,du.
\end{equation}
In particular, if $\inte K=\emptyset$ and $G(0, u)=0$ for $u\in \sphere$, then $\dveV(K)=0$.

The general dual volume was introduced for $K\in \mathcal{K}^n_{(o)}$ in \cite{GHWXY} via (\ref{def-H-volume-06-23}), but there, $\wp$ is a continuous function from $(0, \infty)\times \sphere$ to $(0,\infty)$.  This is clearly subsumed under Definition~\ref{july251}, since such a $G$ can be extended to $[0, \infty)\times \sphere$ by setting $G(0,u)=0$ for $u\in\sphere$.  In fact, the definition in \cite{GHWXY} applies to star-shaped sets about $o$ whose radial functions are positive bounded Borel measurable functions on $\sphere$. A similar extension is possible in Definition~\ref{july251}, but in this paper we only need to work with sets in $\mathcal{K}^n_{o}$.

In Definition~\ref{july251}, one can of course take any continuous function $G(\cdot, \cdot): [0, \infty)\times \sphere\rightarrow [0, \infty)$. For example, when $G(t, u)=t^q/n$ for $q>0$, $$\dveV(K)=\frac{1}{n}\int_{\sphere} \rho_K^q(u)\,du$$
is the {\em $q$th dual volume} of $K$; see \cite[p.~410]{Gar06}. Another example is obtained by setting
$$ G(t, u)=\int^{t}_0\phi(ru)r^{n-1}\,dr$$
for $t\ge 0$ and $u\in \sphere$, where $\phi:\mathbb{R}^n \rightarrow [0, \infty)$ is an integrable function; indeed, in this case $G(t, u)$ is increasing in $t$ and
$$
\dveV(K) =\int_{\sphere}\int^{\rho_{K}(u)}_0\phi(ru)r^{n-1}\,dr\,du= \int_{K}\phi(x)\,dx<\infty.
$$
For more examples, see \cite{GHWXY}.

In \cite[Lemma 6.1]{GHWXY}, $\dveV$ was shown to be continuous on $\mathcal{K}_{(o)}^n$ in the Hausdorff metric when $\wp:(0, \infty)\times \sphere\to(0,\infty)$ is continuous.  We next prove a corresponding result for $\dveV$ on $\mathcal{K}_{o}^n$.  The proof is similar to that of \cite[Lemma~2.2]{BorFor}.

\bl\label{continuity for vG} Let $G: [0, \infty) \times \sphere\rightarrow [0, \infty)$ be continuous.  If $K_i\in \cK_o^n$, $i\in \N$, and $K_i\to K\in\cK_o^n$ as $i\to\infty$, then
$\lim_{i\rightarrow \infty}\dveV(K_{i})=\dveV(K)$.
\el

\begin{proof}
Let $K_i\in \cK_o^n$, $i\in \N$, and $K_i\to K\in\cK_o^n$ as $i\to\infty$.  If $K\in \cK_{(o)}^n$, we can assume without loss of generality that $K_i\in \cK_{(o)}^n$ for all $i\in \N$, and the result then follows from \cite[Lemma~6.1]{GHWXY}.  It therefore suffices to prove the lemma when $o\notin \inte K$.

To this end, suppose first that $\inte K=\emptyset$, so that $K\subset v^{\perp}$ for some $v\in \sphere$. First, we show that $\rho_{K_i}(u)\to \rho_K(u)$ as $i\to\infty$ for $\cH^{n-1}$-almost all $u\in\sphere$. Since $\cH^{n-1}(\sphere\cap v^\perp)=0$, it suffices to consider a fixed $u\in\sphere\setminus v^\perp$.  Let $\varepsilon\in (0,|\langle u,v\rangle|)$.
There exists $i_{\ee}\in\N$ so that $K_i\subset K+\ee^2B^n\subset \{x\in \R^n:\,|\langle x, v\rangle|\leq \ee^2\}$ for $i>i_{\ee}$. Hence, for $i> i_\ee$ we get
$$
0\le \rho_{K_i}(u)\ee< \rho_{K_i}(u)|\langle u,v\rangle|=|\langle \rho_{K_i}(u)u,v\rangle|\le \ee^2,
$$
and therefore $0\le \rho_{K_i}(u)< \ee$ for $i> i_\ee$.  Thus $\rho_{K_i}(u)\to 0=\rho_K(u)$ as $i\to\infty$, as required.

Since $K_i\rightarrow K$, there exists $R>0$ such that $\rho_{K_i}\leq R$ for $i\in \N$.  The continuity of $G$ on $[0, \infty) \times \sphere$ implies that  $M_0=\max\{G(t, u): (t, u)\in [0, R]\times \sphere\}<\infty$.
Hence, since $G$ is continuous, $G(\rho_{K_i}(u),u)\to G(\rho_K(u),u)$ for $\cH^{n-1}$-almost all $u\in\sphere$ and the dominated convergence theorem applies. This yields
$$
\lim_{i\rightarrow \infty} \dveV(K_i)=\lim_{i\to\infty}\int_{\sphere}G(\rho_{K_i}(u),u)\,du=\int_{\sphere}G(\rho_{K}(u),u)\,du=\dveV(K),
$$
as required.

It remains to consider the case when $\inte K\neq \emptyset$ and $o\in \partial K$. It is shown in the proof of \cite[Lemma~2.2]{BorFor} that $\lim_{i\rightarrow \infty}\rho_{K_i}(u)=\rho_K(u)$ for $u\in \sphere\setminus  \partial N(K, o)^*$. Using this and (\ref{dual-cone-6-26}), we obtain
\begin{eqnarray*}
\lim_{i\rightarrow \infty}  \dveV(K_i)   &=&  \lim_{i\rightarrow \infty}   \int_{\sphere\setminus  \partial N(K, o)^*} \wp(\rho_{K_{i}}(u), u)\,du \\
&=&    \int_{\sphere\setminus  \partial N(K, o)^*} \wp(\rho_{K}(u), u)\,du = \dveV(K),
\end{eqnarray*}
where the second equality follows again from the dominated convergence theorem and the fact that $\rho_{K_i}\leq R$ for $i\in \N$.
\end{proof}
	
We shall also need the following lemma.
For  $v\in\sphere$ and $\ee\in (0,1)$, let
\begin{equation}\label{july193}
\Sigma_{\ee}(v)=\{u\in \sphere:\langle u, v\rangle\ge \ee\}.
\end{equation}

\bl \label{decreasing-general-6-30}
Let $\wp(t, u)$ be continuous on $(0, \infty)\times \sphere$ and decreasing in $t$ on $(0, \infty)$. Let $0<\ee_0<1$ and suppose that for $v\in \sphere$,
\begin{equation}\label{condE2-6-30}
\lim_{t\to 0+} \int_{\Sigma_{\ee_0}(v)} \wp(t,u)\,du=\infty.
\end{equation}
If $K_i\in \cK_{(o)}^n$, $i\in \N$, and $K_i\to K\in\cK_o^n$ as $i\to\infty$ with $o\in \partial K$, then
$$\lim_{i\rightarrow \infty} \dveV(K_i)= \infty.$$
\el

\begin{proof}
Let $K_i\in \cK_{(o)}^n$, $i\in \N$, and $K_i\to K\in\cK_o^n$ as $i\to\infty$ with $o\in \partial K$. Choose $v\in \sphere\cap N(K,o)$. Let $t\in (0,\ee_0)$ . Then there is an $i_t\in\N$ such that for $i> i_t$, we have $K_i\subset\{z\in\R^n:\langle z,v\rangle\le t^2\}$. If $u\in\Sigma_{\ee_0}(v)$ and $i> i_t$, then
$$t<\ee_0\le\langle u,v\rangle =\frac{\langle \rho_{K_i}(u)u,v\rangle}{\rho_{K_i}(u)}\le \frac{t^2}{\rho_{K_i}(u)},$$
and therefore $0\le \rho_{K_i}(u)< t$ for  $i> i_t$ and $u\in\Sigma_{\ee_0}(v)$. Hence, for $i> i_t$ we get
\begin{equation}\label{upper-bounded-06-30}
\dveV(K_i)\geq \int_{\Sigma_{\ee_0}(v)} \wp(\rho_{K_{i}}(u), u)\,du\geq \int_{\Sigma_{\ee_0}(v)} \wp(t, u)\,du,
\end{equation}
since $G(\cdot,u)$ is decreasing for $u\in\sphere$. In (\ref{upper-bounded-06-30}), we let $i\rightarrow \infty$ and then $t\rightarrow 0+$.  This and (\ref{condE2-6-30}) yield the assertion.
\end{proof}
	
\section{The general dual Orlicz Minkowski problem for discrete measures}\label{discrete}

Recall that $\pmb{\alpha}^*_K$ was defined by (\ref{july151}) and $\alpha_K$ was also introduced in Section~\ref{Maps}.

Let $K\in \cK_{(o)}^n$ and let $\wp: (0, \infty)\times \sphere\rightarrow (0, \infty)$ and $\psi: (0, \infty)\rightarrow (0, \infty)$ be continuous.  Suppose that $\wp_t(t,u)=\partial \wp(t,u)/\partial t$ is such that $u\mapsto\wp_t(\rho_K(u),u)$ is integrable on $\sphere$. The {\em general dual Orlicz curvature measure} $\deV(K, \cdot)$ is defined (see \cite[Definition~3.1]{GHWXY}) by
\begin{equation}\label{gencdef-06-25}
\deV(K, E)=\frac{1}{n}\int_{\pmb{\alpha}^*_K(E)} \frac{\rho_{K}(u)\, \wp_t(\rho_K(u), u) }{\psi(h_{K}(\alpha_K(u)))}\,du
\end{equation}
for each Borel set $E \subset \sphere$. If $\psi\equiv1$, we often write $\leV(K, \cdot)$ instead of $\deV(K, \cdot)$. It was shown in \cite[Section~3]{GHWXY} that $\deV(K, \cdot)$ is a finite {\em signed} Borel measure on $\sphere$ and several desirable properties were established in \cite[Proposition 6.2]{GHWXY}.

In Definition~\ref{general dual Orlicz curvature measure-7-01}, we shall extend the definition of $\deV(K, \cdot)$ to $K\in \cK_{o}^n$.

Let $g\in C(\Omega)$ and $h_0\in C^+(\Omega)$, where $\Omega\subset\sphere$ is a closed set not contained in any closed hemisphere, and let $a\in \R\cup\{-\infty\}$. Suppose that $\varphi_0\in {\mathcal{J}}_a$ (see (\ref{aug221})) is continuously differentiable and such that $\varphi_0'$ is nonzero on $(0,\infty)$. If $\wp$ and $\wp_t$ are continuous on $(0, \infty)\times \sphere$, then the variational formula from \cite[Theorem 5.3]{GHWXY} asserts that
\begin{equation}\label{variation-11-27-1}
\lim_{\ee\rightarrow 0} \frac{\dveV([h_{\ee}])-\dveV([h_{0}])}{\ee}=
n\int_{\Omega} g(u)\, d\deVo([h_{0}], u),
\end{equation}
where $h_{\ee}$ is given by \eqref{genplus}, the Alexandrov body is taken with respect to $\Omega$, and $\psi_0(t)=t\varphi_0'(t)$.

The following Minkowski-type problem was proposed in \cite[Problem 6.3] {GHWXY}.

\begin{problem}\label{Minkowski-c-06-25}
For which nonzero finite Borel measures $\mu$ on $\sphere$ and continuous functions $\wp:(0, \infty)\times \sphere\rightarrow (0, \infty)$ and $\psi:(0,\infty)\to(0,\infty)$ do there exist $\tau\in \R$ and $K\in \cK_{(o)}^n$ such that $\mu=\tau\,\deV(K,\cdot)$?
\end{problem}

In \cite{GHWXY}, it was shown that this requires finding, for given $G$, $\psi$, and $f:\sphere\to [0,\infty)$, an $h:\sphere\to (0,\infty)$ and $\tau\in \R$ that solve the Monge-Amp\`{e}re equation
\begin{equation}\label{MPnew}
\frac{\tau h}{\psi\circ h}\,P(\bar{\nabla}h+h\iota)
\,\det(\bar{\nabla}^2h+hI)=f,
\end{equation}
where $P(x)=|x|^{1-n}G_t(|x|,\bar{x})$, $\bar{\nabla}$ and $\bar{\nabla}^2$ are the gradient vector and Hessian matrix of $h$, respectively, with respect to an orthonormal frame on $\sphere$, $\iota$ is the identity map on $\sphere$, and $I$ is the identity matrix.

A solution to Problem~\ref{Minkowski-c-06-25} was presented in \cite[Theorem~6.4]{GHWXY}, assuming that $G_t<0$, $G$ satisfies some growth conditions, and
\begin{equation}\label{June-22}
\int_{1}^\infty\frac{\psi(s)}{s}\, ds=\infty.
\end{equation}
Our aim here is to provide a solution to Problem~\ref{Minkowski-c-06-25} when $G_t>0$.  To do this, we first deal with the case when $\mu$ is discrete, i.e., when $\mu=\sum_{i=1}^m \lambda_i \delta_{u_i}$, where $\delta_{u_i}$ denotes the Dirac measure at $u_i\in \sphere$, $\lambda_i>0$ for each $i$, and $\{u_1,\dots,u_m\}$ is not contained in a closed hemisphere.  In this case we seek a solution for which $K\in \cK_{(o)}^n$ is a convex polytope.  This discrete Minkowski-type problem  has been solved in several special cases. Indeed, when $G(t, u)=t^n/n$, then $\dveV(K)$ is the volume of $K$ and the corresponding Orlicz-Minkowski problem for discrete measures was solved in \cite{huang, HugLYZ,  liaijun2014}. When $\psi(t)=t^p$ for $p>1$ and $G(t, u)=t^q\phi_1(u)$ for $q>0$ and $\phi_1\in C^+(\sphere)$, the problem becomes the $L_p$ dual Minkowski problem for discrete measures proposed in \cite{LYZ-Lp} and solved in \cite{BorFor}.  Our solution to
Problem~\ref{Minkowski-c-06-25} in Theorem~\ref{whole-general-solution-7-5} below significantly extends those in \cite{BorFor, huang, HugLYZ,  liaijun2014}.

We utilize the techniques in \cite{HugLYZ} which were also found effective in other works, such as \cite{BorFor, HongYeZhang-2017, huang}.  Let $m>n$ be an integer and suppose that  $\{u_1,\dots,u_m\}$ is not contained in a closed hemisphere. For each $z=(z_1, \dots, z_m)\in [0,\infty)^m$, let
\begin{equation}\label{july172}
P(z) =\{x\in\R^n:\,\langle x, u_i\rangle\leq z_i, \ \mbox{for}\ i=1,\dots, m\}.
\end{equation}
Then $P(z)\in\cK_o^n$ is a convex polytope. We point out that the facets of $P(z)$ are among the support sets  $F(P(z),u_i)$, $i\in \{1,\ldots,m\}$, of $P(z)$,
but not all of these  necessarily are facets.  We  have $h(P(z),u_i)\le z_i$ for $i\in \{1,\ldots,m\}$, with equality if $F(P(z),u_i)$ is a facet of $P(z)$.

\bl \label{facet-06-27}
Let $P\in \cK_{(o)}^n$ be a convex polytope with facets $F(P, u_1),\dots,  F(P, u_m)$.  Then
\begin{equation}\label{polytope-06-25}
\deV(P, \cdot)=\sum_{i=1}^m \gamma_i\delta_{u_i},
\end{equation}
where
\begin{equation}\label{polytope-06-26}
\gamma_i=\frac{\leV(P, \{u_i\})}{\psi(h_{P}(u_i))}
\end{equation}
for $i=1,\dots,m$.  If $\wp_t>0$ (or $\wp_t<0$) on $(0, \infty)\times \sphere$, then $\gamma_i>0$ (or $\gamma_i<0$, respectively) for $i=1,\dots,m$.
\el

\begin{proof}
That $\deV(P, \cdot)$ is of the form (\ref{polytope-06-25}) follows immediately from the absolute continuity of $\deV(P, \cdot)$ with respect to $S(P,\cdot)$ (see \cite[Proposition~6.2(i)]{GHWXY}), since the latter measure is concentrated on $\{u_1,\dots,u_m\}$.
Using (\ref{gencdef-06-25}) and the fact that by (\ref{july151}), $\pmb{\alpha}^*_K(\{u_i\})=\tilde{\pi}(F(P, u_i))$,  we obtain
\begin{eqnarray}\label{july152}
\gamma_i =\deV(P, \{u_i\}) &=&
\frac{1}{n}\int_{\tilde{\pi}(F(P, u_i))} \frac{\rho_{P}(u)\, \wp_t(\rho_P(u), u) }{\psi(h_{P}(u_i))}\,du \nonumber \\&=&\frac{1}{n\psi(h_{P}(u_i))} \int_{\tilde{\pi}(F(P, u_i))} \rho_{P}(u)\, \wp_t(\rho_P(u), u)  \,du =\frac{\leV(P, \{u_i\})}{\psi(h_{P}(u_i))},
\end{eqnarray}
proving (\ref{polytope-06-26}).  (Recall that $\leV(P, \cdot)$ denotes $\deV(P, \cdot)$ with $\psi\equiv1$ and note that $h_P(u_i)>0$, so $\psi(h_P(u_i))>0$ is defined.)

Suppose that $\wp_t>0$ (or $\wp_t<0$) on $(0, \infty)\times \sphere$.  Since $o\in \inte P$, we have $\rho_P>0$, so the integrand in (\ref{july152}) is positive on $\sphere$ (or negative on $\sphere$, respectively).  It follows that $\gamma_i>0$ (or $\gamma_i<0$, respectively) for $i=1,\dots,m$.
\end{proof}

\bl\label{ineqlem}
Let $f: (0, \infty)\rightarrow (0, \infty)$ be continuously differentiable and let $\alpha_2>2\alpha_1>0$.  There exists $c_0=c_0(\alpha_1,\alpha_2)>0$ such that
\begin{equation}\label{rho0half}
f(\alpha-s)\ge f(\alpha)-c_0s
\end{equation}
for $\alpha\in [2\alpha_1, \alpha_2]$ and $s\in [0, \alpha_1]$.
\el

\begin{proof}
If $\alpha\in [2\alpha_1, \alpha_2]$ and $s\in [0, \alpha_1]$, then $\alpha-s\in [\alpha_1, \alpha_2]$.  Let $c_0=\max\{|f'(s)|: s\in [\alpha_1, \alpha_2]\}$.  Define $g(s)=f(\alpha-s)-f(\alpha)+c_0s$ for $s\ge 0$. Then $g(0)=0$ and $g'(s)=c_0-f'(\alpha-s)\ge 0$ for $s\in [0, \alpha_1]$.  Therefore on $[0, \alpha_1]$,  $g$ is increasing and hence $g(s)\ge g(0)=0$, which proves (\ref{rho0half}).
\end{proof}

For $K\in\K^n_{(o)}$, we let
\begin{equation}\label{july176}
\|h_{K}\|_{\varphi, \mu}=\inf\left\{ \lambda>0: \frac{1}{\varphi(1)\,\mu(\sphere)}\int _{\sphere} \varphi\left(\frac{h_{K}(u)}{\lambda}\right)\,d\mu(u) \leq 1\right\}.
\end{equation}
Under appropriate assumptions, $\|\cdot\|_{\varphi, \mu}$ is a norm, the Orlicz or Luxemburg norm.  For example, in \cite[Section~4]{GHW2014}, and elsewehere, the triangle inequality is proved for suitable convex $\varphi$.  We do not need this restriction on $\varphi$, since we merely use (\ref{july176}) for normalization purposes.  Hereafter, we shall use $\|\cdot\|_{\varphi, \mu}$ for a nonzero finite Borel measure $\mu$ on $\sphere$ and a continuous, strictly increasing function $\varphi:[0,\infty)\to [0,\infty)$ with $\varphi(0)=0$ and $\varphi(t)\to\infty$ as $t\to\infty$.

Suppose that $\psi: (0, \infty)\rightarrow (0, \infty)$ is continuous.  Define
\begin{equation}\label{def-varphi-06-26} \varphi(t)=\int_{0}^t\frac{\psi(s)}{s}\, ds ~\quad{\text{for}}~~t>0~~\quad{\text{and}}~~\quad\varphi(0)=0.
\end{equation}
(A similar, but slightly different, function was employed in \cite[(65)]{GHWXY}.) If $\varphi<\infty$ on $(0, \infty)$, then it is continuous (by the dominated convergence theorem) and strictly increasing on $[0,\infty)$, and $\varphi'(t)=\psi(t)/t$ for $t>0$. Note that this assumption on $\varphi$ imposes a weak growth condition on $\psi(t)$ as $t\downarrow 0$.

The hypotheses of the next theorem allow $\psi(t)=t^p$ for $p>1$ and $G(t,u)=t^q$ for $q>0$, for example.

\bt \label{whole-solution-6-27}
Let $\mu=\sum_{i=1}^m \lambda_i \delta_{u_i}$, where $\lambda_i>0$, $i=1,\dots,m$, and $\{u_1,\dots, u_m\}\subset \sphere$ is not contained in a closed hemisphere. Let $\wp: [0, \infty)\times \sphere\rightarrow [0, \infty)$ be continuous and such that  $\wp_t$ is continuous and positive on  $(0, \infty)\times \sphere$.  Suppose that $\psi: (0, \infty)\rightarrow (0, \infty)$ is continuous and such that $\lim_{t\rightarrow 0+} \psi(t)/t=0$ and \eqref{June-22} holds.  Then there exist a convex polytope $P\in \cK_{(o)}^n$ and $\tau>0$ such that
\begin{equation}\label{july282}
\mu=\tau\,\deV(P, \cdot) \ \ \mathrm{and} \ \ \|h_P\|_{\varphi, \mu}=1,
\end{equation}
where $\varphi$ and $\tau$ are given by \eqref{def-varphi-06-26} and \eqref{solution-explicit-6-30}, respectively.
\et

\begin{proof}
Define $\varphi$ by (\ref{def-varphi-06-26}). The assumptions on $\psi$ imply that $\varphi$ is finite, so
$\varphi:[0,\infty) \to [0,\infty)$ is continuous and strictly increasing, and that $\varphi(t)\to\infty$ as $t\to\infty$. It follows that the set
\begin{equation}\label{july161}
M=\left\{(z_1, \dots, z_m)\in [0,\infty)^m: \sum_{i=1}^m\lambda_i \varphi(z_i)=\sum_{i=1}^m\lambda_i  \varphi(1)\right\}
\end{equation}
is compact and nonempty as $(1,\dots, 1)\in M$.  By Lemma~\ref{continuity for vG} and since $z\mapsto P(z)$, $z\in [0,\infty)^m$, is continuous, there is a $z^0=(z_1^0, \dots, z_m^0)\in M$ such that
\begin{equation}\label{optimization-polytope-6-26}
\dveV(P(z^0))= \max\{\dveV(P(z)): \ \ z\in M\}.
\end{equation}
As $\wp_t>0$, $\wp(t, u)$ is strictly increasing in $t\in [0, \infty)$, and then (\ref{def-H-volume-06-23}) implies that $\dveV(\cdot)$ is also increasing, i.e., if $K\subset K'$, then $\dveV(K)\leq \dveV(K')$. From (\ref{july172}) we see that $\ball \subset P((1,\dots, 1))$ and then (\ref{optimization-polytope-6-26}) yields \begin{equation}\label{compare-ball-7-5}
\infty>\dveV(P(z^0))\ge \dveV(P((1, \dots, 1)))\ge \dveV(\ball)=\int_{\sphere}G(1, u)\,du>\int_{\sphere}G(0, u)\,du.
\end{equation}
In view of \eqref{july132}, this implies that $\dim P(z^0)=n$.

We shall first prove the theorem assuming that $o\in \inte P(z^0)$, in which case $P(z^0)\in \cK_{(o)}^n$ and $z_i^0>0$ for $i=1,\ldots,m$. Fix $i\in\{1,\ldots,m\}$ for the moment. Let $h_0\in C^+(\sphere)$ be such that $h_0(u_j)=z_j^0>0$ for $j=1,\ldots,m$. Further, let
$g_i\in C(\sphere)$ be such that $g_i(u_j)=\delta_{ij}$ for $j=1,\ldots,m$. If $|\ee|$ is small enough, we may define $h_{\ee}$ via (\ref{genplus}) with $\varphi_0(t)=t$, so that $h_\ee=h_0+\ee g_i\in C^+(\sphere)$ and $\psi_0(t)=t\varphi_0'(t)=t$. Moreover, if the Alexandrov body $[h_\ee]$
of $h_\ee$ is taken with respect to the set $\Omega=\{u_1,\ldots,u_m\}$, we have $[h_\ee]=P(z^0+\ee e_i)$.

Using this, (\ref{variation-11-27-1}), and (\ref{july152}), we obtain
\begin{eqnarray}\label{variation-11-27-1-1}
\frac{\partial \dveV(P(z))}{\partial z_i}\bigg|_{z=z^0} &=&
\lim_{\ee\rightarrow 0} \frac{\dveV([h_{\ee}])-\dveV([h_{0}])}{\ee} \nonumber \\&=&
n\sum_{j=1}^m g_i(u_j)\deVo(P(z^0), \{u_j\}) \nonumber \\
&=&n\deVo(P(z^0), \{u_i\})=n\frac{\leV(P(z^0), \{u_i\})}{h_{P(z^0)}(u_i)}.
\end{eqnarray}
The argument does not depend on the choice of $z\in\R^m$ with positive coordinates $z_i$, so the calculation shows that the map $z\mapsto \dveV(P(z))$ has continuous partial derivatives and therefore is continuously differentiable. Moreover, since $\lambda_i>0$ and $\varphi'>0$, the rank condition in the Lagrange multiplier theorem is satisfied. In view of (\ref{july161}) and (\ref{optimization-polytope-6-26}), that
theorem provides a $\tau\in \R$ such that
\begin{equation}\label{july171}
\frac{ \tau}{n}\frac{\partial \dveV(P(z))}{\partial z_i}\bigg|_{z=z^0} =\frac{\partial}{\partial z_i}\sum_{i=1}^m\lambda_i \varphi(z_i)\bigg|_{z=z^0}
\end{equation}
for $i=1,\dots,m$.  As $o\in \inte P(z^0)$, (\ref{july172}) implies that $z^0_i>0$ for each $i$.  This, (\ref{variation-11-27-1-1}), (\ref{july171}), and $\varphi'(t)=\psi(t)/t$ for $t>0$ yield
\begin{equation}\label{variation-11-27-1-11}
\tau\frac{\leV(P(z^0), \{u_i\})}{h_{P(z^0)}(u_i)} = \lambda_i \varphi'(z_i^0)=  \lambda_i \frac{\psi(z_i^0)}{z_i^0}~~\quad~~{\text{for}}~~i=1,\dots,m,
\end{equation}
while $z^0\in M$ implies that
\begin{equation}\label{norm-6-28}
\sum_{i=1}^m\lambda_i\varphi(z_i^0)= \varphi(1)\sum_{i=1}^m\lambda_i.
\end{equation}
For each $i$, we have $\lambda_i>0$ and hence $\leV(P(z^0), \{u_i\})>0$, by (\ref{variation-11-27-1-11}).  The absolute continuity of $\deV(P, \cdot)$ with respect to $S(P,\cdot)$ (see \cite[Proposition~6.2(i)]{GHWXY}) implies that the face $F(P(z^0), u_i)$ is actually a facet, hence we have $z^0_i=h_{P(z^0)}(u_i)$.  From (\ref{variation-11-27-1-11}), we conclude that
\begin{equation}\label{july175}
\lambda_i= \tau\frac{\leV(P(z^0), \{u_i\})}{\psi(h_{P(z^0)}(u_i))}= \tau\,\deV(P(z^0), \{u_i\}) \end{equation}
for $i=1,\dots,m$.  This proves that $\mu=\tau\,\deV(P(z^0), \cdot)$ because both measures are concentrated on $\{u_1,\dots,u_m\}$.  Summing (\ref{july175}) over $i$, we obtain
\begin{equation} \label{solution-explicit-6-30}
\tau =\frac{\mu(\sphere)}{\deV(P(z^0),\sphere)}  =\frac{1}{\leV(P(z^0), \sphere)}  \int_{\sphere}\psi(h_{P(z^0)}(u))\,d\mu(u).
\end{equation}
Moreover, in view of (\ref{july176}), (\ref{norm-6-28}) is equivalent to  $\|h_{P(z^0)}\|_{\varphi, \mu}=1$.

This proves the theorem under the assumption that $o\in \inte P(z^0)$, which we now claim is true.  Suppose that $o\in \partial P(z^0)$. To obtain a contradiction, we use an argument similar to that in the proof of \cite[Lemma~3.2]{BorFor}.   By relabeling, if necessary, we may suppose that for some $1\leq k<m$,  $z_j^0=0$ for $j=1,\dots,k$  and $z_j^0>0$  for $j=k+1,\dots, m$. Note that $k<m$ because otherwise, $z_j^0=0$  for $j=1,\dots, m$ implies $P(z^0)=\{o\}$, which is impossible. Let
\begin{equation}\label{july177}
\lambda=\frac{\lambda_1+\cdots +\lambda_k}{\lambda_{k+1}+\cdots+\lambda_m}>0
\end{equation}
and choose $t_0>0$ small enough that $\varphi(z_{i}^0)-\lambda\varphi(t_0)>0$ for $i=k+1, \dots, m$. For $t\in (0, t_0)$, let
\begin{equation}\label{july178}
{a}^t=\left(0,\dots,0,  \varphi^{-1}(\varphi(z_{k+1}^0)-\lambda\varphi(t)), \dots,  \varphi^{-1}(\varphi(z_m^0)-\lambda\varphi(t))\right),
\end{equation}
where the first $k$ components of $a^t$ are equal to $0$, and let \begin{equation}\label{july222}
b^t=a^t+t(e_1+\cdots+e_k),
\end{equation}
so that $b^t$ is obtained from $a^t$ by setting the first $k$ components equal to $t$. By (\ref{def-varphi-06-26}), $\varphi$, and hence $\varphi^{-1}$, is increasing on $[0, \infty)$.  Therefore  $$a^t_i=\varphi^{-1}(\varphi(z_{i}^0)-\lambda\varphi(t))\leq z_i^0$$
for $i=k+1,\dots, m$.  This yields
\begin{equation}\label{inclusion-6-26}
P(a^t)\subset P(z^0)  \ \ \mathrm{and~hence}\ \   \dveV(P(a^t))\leq \dveV(P(z^0)).
\end{equation}
For $t\in (0, t_0)$, we have $o\in \inte P(b^t)$ and from (\ref{july178}) and (\ref{july222}),
\begin{equation}\label{compare-6-26}
P(a^t)\subset P(b^t) \ \ \mathrm{and~hence}\ \  \dveV(P(a^t))\leq \dveV(P(b^t)).
\end{equation}
Using (\ref{july177}), (\ref{july178}), $\varphi(z_1^0)=\cdots=\varphi(z_k^0)=\varphi(0)=0$, and $z^0\in M$, we obtain
\begin{eqnarray*}
\sum_{i=1}^m\lambda_i \varphi(b_i^t) &=&\sum_{i=1}^k \lambda_i  \varphi(t)+\sum_{i=k+1}^m \lambda_i  (\varphi(z_{i}^0)-\lambda\varphi(t)) \\ &=& \varphi(t) \left(\sum_{i=1}^k \lambda_i -\lambda \sum_{i=k+1}^m \lambda_i\right) +\sum_{i=k+1}^m \lambda_i \varphi(z_{i}^0) \\&=&\sum_{i=1}^m \lambda_i\varphi(z_{i}^0) =\sum_{i=1}^m \lambda_i \varphi(1),
\end{eqnarray*}
from which we see via (\ref{july161}) that $b^t\in M$.

Let $r_0=\min\{z_{i}^0: i=k+1,\dots,m\}>0$ and let $R_0> \max\{z_{i}^0: i=k+1,\dots,m\}$ be such that $P(z^0)\subset \inte R_0\ball$.  We apply Lemma~\ref{ineqlem} with $f=\varphi^{-1}$, $\alpha=\varphi(z_i^0)$ for $i=k+1, \dots, m$, $s=\lambda \varphi(t)>0$,  $\alpha_1=\varphi(r_0)/2$, and $\alpha_2=\varphi(R_0)$.  We conclude that with
$$c_0=\max\{(\varphi^{-1})'(s): s\in [\alpha_1, \alpha_2]\} =\max\left\{\frac{1}{\varphi'(\varphi^{-1}(s))}: s\in [\alpha_1, \alpha_2]\right\}$$
and
$$\beta=\min\left\{\varphi^{-1}\left(\frac{\varphi(r_0)}{2\lambda}\right),\, \varphi^{-1}\left(\frac{r_0}{2\lambda c_0}\right)\right\}>0,$$
we have
\begin{equation}\label{bernouli-inequa-6-29}  \varphi^{-1}(\varphi(z_{i}^0)-\lambda\varphi(t)) \geq z_i^0-c_0\lambda\varphi(t)>\frac{r_0}{2}
\end{equation}
for $t\in (0,\beta)$, where the second inequality follows from $z_i^0>r_0$ and the definition of $\beta$.  Note that (\ref{bernouli-inequa-6-29}) and the definition of $t_0$ ensure that we can choose $t_0=\beta$.
	
Recall (see (\ref{aug231})) that $\rho_{P(z^0)}(u)>0$ if and only if $u\in \sphere\cap N(P(z^0),o)^*$.  By (\ref{july172}), the inclusion in (\ref{inclusion-6-26}), and the fact that $a_i^t=0$ if and only if $z_i^0=0$, $\rho_{P(a^t)}(u)>0$ if and only if $\rho_{P(z^0)}(u)>0$.  In fact, by the definition of $r_0$ and $R_0$, (\ref{july178}), and (\ref{bernouli-inequa-6-29}),  we have $\rho_{P(z^0)}(u)\in (r_0, R_0)$ and $\rho_{P(a^t)}(u)\in (r_0/2, R_0)$ for $u\in \sphere\cap N(P(z^0), o)^*$.  Consequently, in view of the continuity of $G_t$ on $(0,\infty)\times \sphere$, there are constants $c_1>0$ and $\beta_1\in (0,\beta)$ such that
\begin{equation}\label{july191}
G(\rho_{P(z^0)}(u), u)-G(\rho_{P(a^t)}(u), u)\leq c_1(\rho_{P(z^0)}(u) -\rho_{P(a^t)}(u))
\end{equation}
for $u\in \sphere\cap N(P(z^0), o)^*$.

Let $u\in \sphere\cap N(P(z^0), o)^*$ and choose $i_0\in\{k+1,\dots, m\}$ so that the ray from $o$ in the direction $u$ meets the facet $F(P(a^t),u_{i_0})$, and hence, by the inclusion in (\ref{inclusion-6-26}), the facet $F(P(z^0),u_{i_0})$ as well.  Then, by (\ref{bernouli-inequa-6-29}),
$$R_0\langle u,u_{i_0}\rangle>\rho_{P(a^t)}(u)\langle u,u_{i_0}\rangle
=a^t_{i_0}=\varphi^{-1}(\varphi(z^0_{i_0})-\lambda \varphi(t))>\frac{r_0}{2}$$
and $z^0_{i_0}=\rho_{P(z^0)}(u)\langle u,u_{i_0}\rangle$.  Using these relations and (\ref{bernouli-inequa-6-29}) again, we obtain
\begin{equation}\label{estimate-radial-difference-6-29}
\rho_{P(a^t)}(u)=\frac{\varphi^{-1}(\varphi(z^0_{i_0})-\lambda \varphi(t))}{\langle u,u_{i_0}\rangle} \geq \frac{ z_{i_0}^0-c_0\lambda\varphi(t)}{\langle u,u_{i_0}\rangle}\geq \rho_{P(z^0)}(u)-\frac{2R_0c_0\lambda}{r_0}\varphi(t).
\end{equation}
From (\ref{july191}) and (\ref{estimate-radial-difference-6-29}) we get
\begin{eqnarray}\label{difference-generaldual-6-30}
\dveV(P(z^0))-\dveV(P(a^t))&=&\int_{\sphere\cap N(P(z^0), o)^*}\left(G(\rho_{P(z^0)}(u), u)-G(\rho_{P(a^t)}(u), u)\right)\,du \nonumber \\&\leq& c_1\int_{\sphere\cap N(P(z^0), o)^*}\left(\rho_{P(z^0)}(u) -\rho_{P(a^t)}(u)\right)\,du \nonumber \\&\leq&
\frac{2R_0c_0c_1\lambda}{r_0}\varphi(t) \int_{\sphere\cap N(P(z^0), o)^*}\,du \le c_2\varphi(t)
\end{eqnarray}
for $t\in (0,\beta_1)$, where $c_2=2R_0c_0c_1\lambda n\kappa_n/r_0$.

Using $o\in\inte P(b^t)$ and the containments in (\ref{inclusion-6-26}) and (\ref{compare-6-26}), one can show that there is a closed set $E_t\subset \sphere$ and constants $r_1>0$,  $\beta_2\in (0, \beta_1)$, and $c_3>0$, depending only on $n$, $r_0$, and $R_0$, satisfying ${\mathcal{H}}^{n-1}(E_t)\ge c_3t$ for $t\in (0,\beta_2)$ and such that $\rho_{P(a^t)}(u)=0$ and $\rho_{P(b^t)}(u)\geq r_1$ for $u\in E_t$. We omit the details, since these are given in the proof of \cite[p.~13]{BorFor}; there, the set $E_t$ is denoted by $\widetilde{G}_t$ and is the radial projection on $\sphere$ of a certain $(n-1)$-dimensional spherical cylinder of height $t$.  For $u\in E_t$, we have $G(\rho_{P(a^t)}(u),u)=G(0,u)$ and $G(\rho_{P(b^t)}(u),u)\ge G(r_1,u)$ as $G_t>0$.
Consequently,
\begin{eqnarray}\label{july192}
\dveV(P(b^t))&=&\int_{\sphere\setminus E_t} G(\rho_{P(b^t)}(u),u) \,du+\int_{E_t} G(\rho_{P(b^t)}(u),u)\,du\nonumber \\ &\geq & \dveV(P(a^t))  +\int_{E_t}\left(G(r_1,u)-G(0, u)\right)\,du\geq \dveV (P(a^t))+c_4t
\end{eqnarray}
for $t\in (0, \beta_2)$, where
$$c_4=c_3\min\{G(r_1,u)-G(0, u):u\in \sphere\}>0.$$
From (\ref{difference-generaldual-6-30}) and (\ref{july192}), we obtain
\begin{equation}\label{difference-lim-6-30}
\liminf\limits_{t\rightarrow 0+} \frac{\dveV(P(b^t))- \dveV(P(z^0))}{t} \geq   \lim_{t\rightarrow 0+}  \frac{c_4 t -c_2\,\varphi(t)}{t} =c_4 >0,
\end{equation}
since
$$\lim_{t\rightarrow 0+} \frac{\varphi(t)}{t}=\lim_{t\rightarrow 0+} \varphi'(t)=\lim_{t\rightarrow 0+} \frac{\psi(t)}{t}=0.$$
By (\ref{difference-lim-6-30}), there exists $t_1\in (0, \beta_2)$ such that $\dveV(P(b^{t_1}))>\dveV(P(z^0))$.  It was shown above that $b^{t_1}\in M$, so this contradicts (\ref{optimization-polytope-6-26}).  Thus $o\in \inte P(z^0)$ and the proof is complete.
\end{proof}

Recall that for $v\in\sphere$ and $\varepsilon\in (0,1)$, $\Sigma_{\varepsilon}(v)$ is defined by (\ref{july193}) and that $\|\cdot\|_{\varphi,\mu}$ is defined by (\ref{july176}).  The hypotheses of the next theorem allow $\psi(t)=t^p$ for $p>0$ and $G(t,u)=t^q$ for $q<0$, for example.

\bt \label{solution-general-dual-Orlicz-main theorem-11-27}
Let $\mu=\sum_{i=1}^m \lambda_i \delta_{u_i}$, where $\lambda_i>0$, $i=1,\dots,m$, and $\{u_1,\dots, u_m\}\subset \sphere$ is not contained in a closed hemisphere.  Let $\wp: (0, \infty)\times \sphere\rightarrow (0, \infty)$ be continuous and such that $\wp_t$ is continuous and negative on $(0, \infty)\times \sphere$.  Let $0<\ee_0<1$ and suppose that \eqref{condE2-6-30} holds for $v\in \sphere$.  Suppose that $\psi: (0, \infty)\rightarrow (0, \infty)$ is continuous, \eqref{June-22} holds, and that $\varphi$ is finite when defined by \eqref{def-varphi-06-26}.  Then there exist a convex polytope $P\in \cK_{(o)}^n$ and $\tau<0$ such that
$$\mu=\tau\,\deV(P, \cdot) \ \ \mathrm{and} \ \ \|h_P\|_{\varphi, \mu}=1,$$
where $\tau$ is given by \eqref{solution-explicit-6-30}.
\et

\begin{proof}
Define $\varphi$ by (\ref{def-varphi-06-26}). The assumption that $\varphi$ is finite implies that $\varphi:[0,\infty) \to [0,\infty)$ is continuous and strictly increasing, and from \eqref{June-22}, we have $\varphi(t)\to\infty$ as $t\to\infty$. The set $$M'=\left\{(z_1, \dots, z_m)\in (0,\infty)^m:   \sum_{i=1}^m\lambda_i \varphi(z_i)=\sum_{i=1}^m\lambda_i  \varphi(1)\right\}$$
is bounded and nonempty as $(1,\dots,1)\in M'$. Let
\begin{equation}\label{optimization-polytope-6-30}
\alpha=\inf \{\dveV(P(z)): \ \ z\in M'\}, \end{equation}
where $P(z)$ is defined by (\ref{july172}).  Choose $z^j \in M'$, $j\in \N$, such that \begin{equation}\label{optimization-limit-6-30}
\lim_{j\rightarrow \infty} \dveV(P(z^j)) =\alpha.
\end{equation}

Since $M'$ is bounded, we can assume, by taking a subsequence, if necessary, that $z^j\rightarrow z^0\in M$, where $M$ is defined by (\ref{july161}).  However, we actually have $o\in \inte P(z^0)$ and hence $z^0\in M'$.  To see this, suppose to the contrary that
$o\in \partial  P(z^0)$.  Since $P(z^j)\in \cK_{(o)}^n$ for $j\in \N$ and $P(z^j)\rightarrow P(z^0)$ as $j\rightarrow \infty$,    Lemma~\ref{decreasing-general-6-30} yields
$$\lim_{j\rightarrow \infty} \dveV(P(z^j))= \infty.$$
By (\ref{july172}), $\ball \subset P((1,\dots, 1))$.  Also, as $G_t<0$, $G(t, \cdot)$ is decreasing on $(0, \infty)$ and hence $\dveV(\cdot)$ is also decreasing, i.e., if $K\subset K'$, then $\dveV(K)\geq \dveV(K')$.  Therefore, using $(1,\dots, 1) \in M'$, (\ref{optimization-polytope-6-30}), and (\ref{optimization-limit-6-30}), we obtain
$$
\infty> \dveV(\ball)\geq \dveV(P((1,\dots, 1)))\geq \alpha=\lim_{j\rightarrow \infty} \dveV(P(z^j)),
$$
a contradiction proving that $z^0\in M'$ and $P(z^0)\in \mathcal{K}_{(o)}^n$.  By \cite[Lemma~6.1]{GHWXY}, $\dveV(\cdot)$ is continuous in the Hausdorff metric on $\mathcal{K}_{(o)}^n$, so
\begin{equation}\label{compare-ball-7-9}
\infty> \dveV(\ball)\ge \dveV(P(z^0))=\lim_{j\rightarrow \infty} \dveV(P(z^j))=\alpha>0.
\end{equation}

The remainder of the proof is precisely the same as the passage from (\ref{variation-11-27-1-1}) to (\ref{solution-explicit-6-30}) in the proof of Theorem~\ref{whole-solution-6-27}.
\end{proof}

Under the conditions on $\mu$, $\wp$, and $\psi$ stated in Theorem~\ref{solution-general-dual-Orlicz-main theorem-11-27}, but with the assumption that $\varphi<\infty$ replaced by the condition
\begin{equation}\label{july195}
\lim_{t\to \infty} \int_{\sphere} \wp(t,u)\,du=0,
\end{equation}
\cite[Theorem~6.4]{GHWXY} proves that there is a $K\in \cK_{(o)}^n$ such that
\begin{equation}\label{msol}
\frac{\mu}{ \mu(\sphere) }=\frac{\deV(K, \cdot)}{\deV(K, \sphere)}.
\end{equation}
If $\mu$ is discrete, \cite[Theorem~6.4]{GHWXY} does not prove that $K$ is a convex polytope, but, as is explained in the discussion after \cite[Corollary~6.5]{GHWXY}, this is an easy consequence of (\ref{msol}).  Thus Theorem~\ref{solution-general-dual-Orlicz-main theorem-11-27} is a variant of \cite[Theorem~6.4]{GHWXY} for discrete $\mu$.

\section{General dual Orlicz curvature measures for compact convex sets}\label{section-5}
The general dual Orlicz curvature measure $\deV(K, \cdot)$ was defined by (\ref{gencdef-06-25}) for $K\in \cK_{(o)}^n$.  In this section, we extend the definition to $K\in \cK_{o}^n$.

Let $K\in\cK_o^n$. Recall that  $N(K,o)$, $N(K,o)^*$, and $\pmb{\alpha}^*_K$ are defined by   (\ref{july201}),  (\ref{july135}), and (\ref{july151}), respectively.

\bd\label{general dual Orlicz curvature measure-7-01}
Define the {\em general Orlicz curvature measure} $\deV(K, \cdot)$ by
\begin{equation}\label{gencdef-7-1-5}
\deV(K, E)=\frac{1}{n}\int_{\pmb{\alpha}^*_K(E\setminus N(K, o))} \frac{\rho_K(u)\,\wp_t(\rho_K(u), u) }{\psi(h_{K}(\alpha_K(u)))}\,du
\end{equation}
for each Borel set $E\subset\sphere$, whenever $\wp: [0, \infty)\times \sphere\rightarrow [0, \infty)$ and $\psi: [0, \infty)\rightarrow [0, \infty)$ with $\psi(t)>0$ for $t>0$ are such that the integral in \eqref{gencdef-7-1-5} exists for all $K\in \mathcal{K}^n_{o}$ and Borel sets $E\subset\sphere$.
\ed

Note that if $\dim K<n$, then $\deV(K, E)=0$, since $\sphere\cap N(K,o)^*$ is then at most $(n-2)$-dimensional. Furthermore, if $\dim K=n$, then in view of (\ref{dual-cone-6-26}) and (\ref{aug232}), the integral in (\ref{gencdef-7-1-5}) may equivalently be taken over $\pmb{\alpha}_K^*(E)\cap N(K,o)^*$ or over $\pmb{\alpha}_K^*(E)\cap \inte N(K,o)^*$. For $\cH^{n-1}$-almost all $u\in \pmb{\alpha}^*_K(E\setminus N(K, o))$, the vector $\alpha_K(u)$ is well defined and $\alpha_K(u)\notin N(K,o)$, hence $h_{K}(\alpha_K(u))>0$ and $\psi(h_{K}(\alpha_K(u)))>0$.

As before, if $\psi\equiv 1$, we often write $\leV(K, \cdot)$ instead of $\deV(K, \cdot)$.
The integral in (\ref{gencdef-7-1-5}) should be considered as $0$ if $\pmb{\alpha}^*_K(E\setminus N(K, o))=\emptyset$, in particular for $E\subset \sphere\cap N(K, o)$.  In other words, $\deV(K, E)=0$ for each Borel set $E\subset \sphere\cap N(K, o)$.

When $o\in \inte K$, we have $N(K, o)=\{o\}$ and hence $E\setminus N(K, o)=E$, so (\ref{gencdef-7-1-5}) agrees with (\ref{gencdef-06-25}).  Moreover, if $\wp: (0, \infty)\times \sphere\rightarrow (0, \infty)$ and $\psi: (0, \infty)\rightarrow (0, \infty)$ are continuous and $\wp_t(t,u)=\partial \wp(t,u)/\partial t$ is such that $u\mapsto\wp_t(\rho_K(u),u)$ is integrable on $\sphere$, then we can extend $\wp$ and $\psi$ by setting $\wp(0,u)=0$ for $u\in\sphere$ and $\psi(0)=0$ and the integral in (\ref{gencdef-7-1-5}) will exist for all $K\in \mathcal{K}^n_{(o)}$ and Borel sets $E\subset \sphere$.  Thus the definition of $\deV(K,\cdot)$ in Section~\ref{discrete} and \cite{GHWXY} is subsumed under Definition~\ref{general dual Orlicz curvature measure-7-01}.

Suppose that $G$ and $\psi$ are such that $\deV(K, \cdot)$ is indeed a finite signed Borel measure on $\sphere$.  Then integrals with respect to $\deV(K, \cdot)$ can be calculated as follows. For any bounded Borel function $g:\sphere\to \R$, we have
\begin{eqnarray}
\int_{\sphere} g(u)\, d\deV(K, u)&=&\frac{1}{n}\int_{\sphere\cap \inte N(K,o)^*}g(\alpha_K(u))\,\frac{\rho_K(u)\,\wp_t(\rho_K(u),u)}
{\psi(h_K(\alpha_K(u)))}\,du \label{new measue-11-27}\\
&=&\frac{1}{n}\int_{\partial K\backslash \Xi_K}g(\nu_{K}(x))\frac{\langle x, \nu_{K}(x)\rangle}{\psi(\langle x, \nu_{K}(x)\rangle)}\,|x|^{1-n}\wp_t(|x|, \tilde{\pi}(x))\,dx.\label{july271}
\end{eqnarray}

Indeed, it suffices to prove (\ref{new measue-11-27}) for $g=1_{E}$, where $E\subset \sphere$ is a Borel set. If $\dim K\le n-1$, then all integrals are zero, so we can assume that $\dim K=n$. Then, using (\ref{relation-11-27}) and (\ref{aug232}), we obtain
\begin{eqnarray*}
\int_{\sphere} 1_E(u)\, d\deV(K, u)&=&\deV(K, E)\\&=& \frac{1}{n}\int_{\pmb{\alpha}_K^*(E)\cap\inte N(K,o)^*} \frac{\rho_K(u)\,\wp_t(\rho_K(u),u)}
{\psi(h_K(\alpha_K(u)))}\,du \\&=&\frac{1}{n}\int_{\sphere\cap \inte N(K,o)^*}1_E(\alpha_K(u))\,\frac{\rho_K(u)\,\wp_t(\rho_K(u),u)}
{\psi(h_K(\alpha_K(u)))} \,du,
\end{eqnarray*}
as required, thus proving (\ref{new measue-11-27}).  Now (\ref{july252}) and a standard change of variables (see, e.g., \cite[(21)]{BorFor}, \cite[(23)]{GHWXY}, or \cite[(2.30)]{LYZActa}) gives (\ref{july271}).  When $\psi$ does not vanish or when $\cH^{n-1}(\Xi_K)=0$, (\ref{july271}) becomes
\begin{equation}\label{form-11-29}
\int_{\sphere} g(u)\, d\deV(K, u)=\frac{1}{n}\int_{\partial K}g(\nu_{K}(x))\frac{\langle x, \nu_{K}(x)\rangle}{\psi(\langle x, \nu_{K}(x)\rangle)}\,|x|^{1-n}\wp_t(|x|, \tilde{\pi}(x))\,dx,
\end{equation}
since it is easy to see that $\langle \nu_K(x),x\rangle=0$ for $x\in\Xi_K\cap \partial K$. We emphasize that these observations are made under the assumption that the integrals exist.

In Proposition~\ref{property-extended-7-1}(iii) below, we will find use for a simplified version of (\ref{new measue-11-27}) that holds when $\psi\equiv 1$ and $tG_t(t, u)=0$ at $t=0$ for $u\in \sphere$.  The latter of these two conditions simply means that $tG_t(t, u)\to 0$ as $t\to 0+$ for $u\in \sphere$.  The apparently weaker condition that $\lim_{t\to 0+}tG_t(t, u)$ exists for $u\in \sphere$ is in fact equivalent.  Indeed, suppose that $u\in\sphere$ and $tG_t(t,u)\to c\neq 0$ as $t\to 0+$. If $c>0$, there exist $0<c_1\le c_2<\infty $ and $t_0>0$ such that $0<c_1\le tG_t(t,u)\le c_2$ for $t\in (0,t_0]$. If $s\in (0,t_0]$, we can divide by $t$ and integrate from $s$ to $t_0$ to obtain
$$G(t_0,u)-c_2\ln t_0+c_2\ln s\le G(s,u)\le G(t_0,u)-c_1\ln t_0+c_1\ln s.$$
But then $G(s,u)\to -\infty$ as $s\to 0+$, a contradiction. If $c<0$, there exist $c_1\le c_2<0 $ and $t_0$ as above and a similar argument leads to $G(s,u)\to \infty$ as $s\to 0+$, again a contradiction.

The following proposition focuses on the case when $\psi\equiv 1$.  In this case, provided $tG_t(t, u)=0$ at $t=0$ for $u\in \sphere$, (\ref{new measue-11-27}) simplifies to
\begin{equation}\label{sept81}
\int_{\sphere} g(u)\, d\leV(K, u)=\frac{1}{n}\int_{\sphere}g(\alpha_K(u))\,\rho_K(u)\,\wp_t(\rho_K(u),u)
\,du
\end{equation}
for any bounded Borel function $g:\sphere\to \R$, since the integral may be restricted to $\inte N(K,o)^*$ due to (\ref{dual-cone-6-26}) and (\ref{july221}).

\bp\label{property-extended-7-1}
Let $\wp: [0, \infty)\times \sphere\rightarrow [0, \infty)$ and let $K\in \cK_o^n$. The following statements hold.

\noindent{\rm{(i)}}  $\leV(K, \cdot)$ is a finite signed measure on $\sphere$.

\noindent{\rm{(ii)}} Suppose that $t^{1-n}\wp_t(t,u)$ is continuous on $[0, \infty)\times \sphere$, where the value of $t^{1-n}\wp_t(t,u)$ for each $u\in \sphere$ at $t=0$ is taken to be the value of the limit as $t\rightarrow 0+$. Then $\leV(K, \cdot)$ is absolutely continuous with respect to $S(K,\cdot)$

\noindent{\rm{(iii)}} Suppose that $tG_t(t,u)$ is continuous on $[0, \infty)\times \sphere$, where $tG_t(t, u)=0$ at $t=0$ for $u\in \sphere$.  If $K_i\in \cK_o^n$ and $K_i\to K\in \cK_o^n$ as $i\to\infty$, then $\leV(K_i,\cdot)\rightarrow \leV(K,\cdot)$ weakly as $i\to\infty$.
\ep

\begin{proof}
(i) As was pointed out before, $\leV(K, \cdot)$ is the zero measure if $\dim K\le n-1$. Hence let $\dim K=n$.  The assumption on $\wp_t$ ensures that with $\psi\equiv 1$, the integral in (\ref{gencdef-7-1-5}) exists.  For the $\sigma$-additivity, it  suffices to show that
\begin{equation}\label{july253}
\leV(K,\cup_{i=1}^{\infty} E_i) =\sum_{i=1}^\infty\leV(K,E_i),
\end{equation}
for disjoint Borel sets $E_i \subset \sphere$, $i\in \N$.  To this end, observe that it follows directly from (\ref{july151}) that $\pmb{\alpha}^*_K((\cup_{i=1}^{\infty}E_i)\setminus N(K, o))=\cup_{i=1}^{\infty} \pmb{\alpha}^*_K(E_i \setminus N(K, o))$.  From (\ref{gencdef-7-1-5}), we see that (\ref{july253}) will be proved if we can show that
\begin{equation}\label{sept82}
{\cH}^{n-1}\big(\pmb{\alpha}^*_K(E_i \setminus N(K, o))\cap \pmb{\alpha}^*_K(E_j \setminus N(K, o))\big)=0
\end{equation}
for $i\neq j$. To see this, note that since $r_K$ is locally bi-Lipschitz on $\sphere\cap \inte N(K,o)^*$, we have
$${\cH}^{n-1}(r_K^{-1}(\partial K\setminus \reg K)\cap \inte N(K,o)^*)=0.$$
Using this and (\ref{aug232}), we get
\begin{eqnarray*}
\lefteqn{{\cH}^{n-1}\big(\pmb{\alpha}^*_K(E_i \setminus N(K, o))\cap \pmb{\alpha}^*_K(E_j \setminus N(K, o))\big)}\\
&=&{\cH}^{n-1}\big(\pmb{\alpha}^*_K(E_i)\cap \pmb{\alpha}^*_K(E_j) \cap \inte N(K,o)^*\cap r_K^{-1}(\reg K) \big).
\end{eqnarray*}
But the latter set is empty, because if it contained a point $u$, we would have
$$r_K(u)\in \pmb{\nu}_K^{-1}(E_i)\cap \pmb{\nu}_K^{-1}(E_j)\cap \reg K=\emptyset,$$
as $E_i\cap E_j=\emptyset$.  This proves (\ref{sept82}) and hence (\ref{july253}).

\smallskip

(ii) If $\inte K=\emptyset$, then $\leV(K, \cdot)=0$ and there is nothing to prove.  Suppose that $\inte K\neq \emptyset$.  Let $E\subset \sphere$ be a Borel set such that $S(K, E)=0$, let $g=1_E$, and choose $R<\infty$ such that $K\subset R\ball$. By (\ref{form-11-29}) with $\psi\equiv 1$, the continuity of $t^{1-n}\wp_t(t,u)$, and the fact that $\langle x, \nu_K(x)\rangle\leq R$ for $x\in\partial K$, we obtain
\begin{eqnarray*}
\leV(K, E) &=& \frac{1}{n}\int_{\partial K}1_E(\nu_{K}(x))\,\langle x, \nu_{K}(x)\rangle\,|x|^{1-n}\wp_t(|x|, \tilde{\pi}(x))\,dx\\ &\leq& \frac{R}{n} \max\left\{t^{1-n}\wp_t(t,u)
: (t, u)\in [0, R]\times \sphere\right\} \cH^{n-1}\left(\{x\in\partial K: \nu_K(x)\in E\}\right)\\
&=& 0,
\end{eqnarray*}
as required.
		
\smallskip

(iii)  The case when $K\in \cK_{(o)}^n$ was proved in \cite[Proposition~6.2(ii)]{GHWXY}.  First assume that $o\in \partial K$ and $\inte K\neq \emptyset$.  Let $g\in C(\sphere)$ and let
$$I_K(u)=g(\alpha_K(u))\,\rho_{K}(u)\, \wp_t(\rho_K(u),u)$$
be the integrand of the right-hand side of (\ref{sept81}). If $u\in \inte N(K,o)^*$, then $u\in \inte N(K_i,o)^*$
for $i\ge i_u$ and $\rho_{K_i}(u)\to \rho_K(u)$ as $i\to\infty$. Let $Z$ be the set consisting of those $u\in\sphere\cap\inte N(K,o)^*$ for which $\rho_K(u)u\notin\reg K$ and those $u\in\sphere\cap\inte N(K_i,o)^*$ for which $\rho_{K_i}(u)u\notin\reg K_i$ for some $i\in\N$. Then (\ref{zerosets}) yields $\cH^{n-1}(Z)=0$.  Also, since $\alpha_{K_i}(u)\to\alpha_K(u)$ as $i\to\infty$ for $u\in \inte N(K,o)^*\setminus Z$ (cf.~\cite[Lemma~2.2]{LYZActa}), we have $I_{K_i}(u)\to I_K(u)$ as $i\to\infty$ for $u\in \inte N(K,o)^*\setminus Z$.

On the other hand, if $u\in\sphere\setminus N(K,o)^*$, then $\rho_K(u)=0$ by (\ref{aug231}) and $\rho_{K_i}(u)\to 0$ as $i\to\infty$ (as can be seen by a separation argument), and hence, using the assumption that $tG_t(t, u)=0$ at $t=0$ for $u\in \sphere$, we have $I_{K_i}(u)\to 0$ as $i\to\infty$.  Thus we have shown that $I_{K_i}(u)\to I_K(u)$ as $i\to\infty$ for $\cH^{n-1}$-almost all $u\in\sphere$.  We also have $\sup\{|I_K(u)|: u\in \sphere\}|<\infty$, by the continuity of $tG_t(t,u)$ on $[0, \infty)\times \sphere$.

Using these facts, (\ref{sept81}), and the dominated convergence theorem, we obtain
\begin{eqnarray*}
\int_{\sphere} g(u)\, d\leV(K, u)&=&\frac{1}{n}\int_{\sphere}I_K(u)\,du\\
&=& \lim_{i\to\infty}\frac{1}{n}\int_{\sphere} I_{K_i}(u)\,du=\lim_{i\to\infty}\int_{\sphere} g(u)\, d\leV(K_i, u),
\end{eqnarray*}
proving the result when $\inte K\neq \emptyset$.

Now assume that $\inte K=\emptyset$.  Since $g$ is continuous on $\sphere$, it is bounded, so there is a $c>0$ such that $|g(\alpha_{K_i}(u))|\le c$ for $u\in\sphere$ and $i\in \N$. We apply Lemma~\ref{continuity for vG} with $G(t,u)$ replaced by $t|\wp_t(t,u)|$ to obtain
\begin{eqnarray*}
\limsup_{i\to\infty}\int_{\sphere} |g(\alpha_{K_i}(u))|\,\rho_{K_i}(u)\,|\wp_t(\rho_{K_i}(u),u)|\,du&\le &
c\limsup_{i\to\infty}\int_{\sphere} \rho_{K_i}(u)\,|\wp_t(\rho_{K_i}(u),u)|\,du\\
&\le &
c\int_{\sphere} \rho_{K}(u)\,|\wp_t(\rho_{K}(u),u)|\,du=0.
\end{eqnarray*}
where we have used again the assumption that $tG_t(t, u)=0$ at $t=0$ for $u\in \sphere$.  This and (\ref{sept81}) yield
$$
\lim_{i\to\infty}\int_{\sphere} g(u)\, d\leV(K_i, u)=0=\int_{\sphere} g(u)\, d\leV(K, u),
$$
completing the proof.
\end{proof}

Finally, we provide a generalization of uniqueness results for $\leV(K, \cdot)$ in \cite[Theorem~6.1]{XY2017-1}, \cite[Theorem~5.2]{zhao}, and \cite[Theorem~3.1]{ZSY2017}, with a simpler proof. We start with a variant of \cite[Lemma~5.1]{zhao}, whose proof we omit since it is similar.

\bl \label{lemma5.3new}
Let $K_1,K_2\in\cK_{(o)}^n$ and let $E=\{v\in\sphere: h_{K_1}(v)>h_{K_2}(v)\}\neq\emptyset$.

\noindent{\rm{(i)}} If $u\in\pmb{\alpha}^*_{K_1}(E)$, then $\rho_{K_1}(u)>\rho_{K_2}(u)$;

\noindent{\rm{(ii)}} $\pmb{\alpha}^*_{K_1}(E)\subset \pmb{\alpha}^*_{K_2}(E)$;

\noindent{\rm{(iii)}} $\mathcal{H}^{n-1}(\pmb{\alpha}^*_{K_2}(E))>0$.
\el

\bt \label{theorem5.4new}
Let $G:(0,\infty)\times\sphere\to (0,\infty)$ be continuous and such that $G_t$ is continuous and negative on $(0,\infty)\times\sphere$.
Suppose that $tG_t(t,u)$ is strictly increasing in $t$ for $u\in\sphere$. If $K,L\in\cK_{(o)}^n$ satisfy $\leV(K, \cdot)=\leV(L, \cdot)$, then $K=L$.
\et

\begin{proof}
Suppose that $K\neq L$.  Then we may assume that $L\not\subset K$.  Let $E=\{v\in\sphere: h_{L}(v)>h_{K}(v)\}\neq\emptyset$. We apply Lemma~\ref{lemma5.3new} with $K_1=L$ and $K_2=K$. Using \eqref{sept81} and the fact that $G_t<0$, Lemma~\ref{lemma5.3new}(i), the assumption that $tG_t(t,u)$ is strictly increasing in $t$ for $u\in\sphere$, and Lemma~\ref{lemma5.3new}(ii) together with \eqref{relation-11-27}, we obtain
\begin{align*}
\leV(K,E)&=\leV(L,E)\\
&=\int_{\sphere} 1_E(\alpha_L(u))\,\rho_L(u)\,G_t(\rho_L(u),u)\, du\\
&\ge \int_{\sphere} 1_E(\alpha_L(u))\,\rho_K(u)\,G_t(\rho_K(u),u)\, du\\
&\ge \int_{\sphere} 1_{\pmb{\alpha}^*_K(E)}(u)\,\rho_K(u)\,G_t(\rho_K(u),u)\, du\\
&=\leV(K, E).
\end{align*}
If $\mathcal{H}^{n-1}(\pmb{\alpha}^*_K(E)\setminus\pmb{\alpha}^*_L(E))>0$, then the second inequality is strict. If $\mathcal{H}^{n-1}(\pmb{\alpha}^*_K(E)\setminus \pmb{\alpha}^*_L(E))=0$,
then Lemma~\ref{lemma5.3new}(iii) implies that $\mathcal{H}^{n-1}(\pmb{\alpha}^*_L(E))>0$ and therefore the first
inequality is strict. Thus, in any case we arrive at a contradiction.
\end{proof}

\section{Minkowski problems for general measures}\label{section-6}
In view of (\ref{gencdef-7-1-5}), one sees that
$$\frac{\,d\deV(K, u)}{\,d\leV(K, u)}=\frac{1}{\psi(h_K(u))}.$$
We consider the following Minkowski-type problem.

\begin{problem}\label{Minkowski-c-07-5}
For which nonzero finite Borel measures $\mu$ on $\sphere$ and continuous functions $\wp: (0, \infty)\times \sphere\rightarrow (0, \infty)$ and $\psi:(0,\infty)\to(0,\infty)$ do there exist $\tau\in \R$ and  $K\in \cK_{o}^n$ with $\inte K\neq \emptyset$ such that $$\mu=\tau\,\deV(K,\cdot) \ \ \mathrm{and/or}\ \ (\psi\circ h_{K})\mu=\tau\,\leV(K,\cdot)?$$
\end{problem}

For our contribution to this problem, we need the following lemma.  It is essentially known (see e.g., \cite{HongYeZhang-2017, huang, LZH2017-a, ZHY2016}), but we provide an explicit dependence of $R$ on $\mu$ that will be needed in the proof of Theorems~\ref{whole-general-solution-7-5} and~\ref{solution-general-dual-Orlicz-main theorem-11-27-even}.

\bl \label{august24}
Let $\mu$ be a finite Borel measure on $\sphere$ not concentrated on any closed hemisphere and let $\varphi:[0,\infty)\to [0,\infty)$ with $\varphi(0)=0$ be continuous and strictly increasing.  Suppose that $K\in \cK_{(o)}^n$ satisfies $\|h_{K}\|_{\varphi, \mu}=1$ (see \eqref{july176}).  Then there is an $R=R(\mu,\varphi)>0$ such that $K\subset RB^n$.
\el

\begin{proof}
There is a $\delta=\delta(\mu)>0$ such that
\begin{equation}\label{posa}
\int_{\sphere}\langle u,v\rangle_+\, d\mu(u)\ge \delta\, \mu(\sphere)
\end{equation}
for $v\in\sphere$, since the integral on the left is continuous in $v$ on $\sphere$ and $\mu$ is not concentrated on a closed hemisphere.  Let $K\in \cK_{(o)}^n$ satisfy $\|h_{K}\|_{\varphi, \mu}=1$ and let $rv_0\in K$, where $r\ge 0$ and $v_0\in\sphere$. Then $[o,rv_0]\subset K$ implies that $h_K(u)\ge r\langle u,v_0\rangle_+$ for $u\in\sphere$, so using \eqref{july193}, \eqref{july176}, and our assumptions on $\varphi$, we obtain
\begin{align}\label{posb}
\varphi(1)\,\mu(\sphere)=  \int _{\sphere} \varphi(h_{K}(u))\,d\mu(u)
\ge  \int _{\Sigma_{\delta/2}(v_0)} \varphi(r\delta/2)\,d\mu(u)
= \varphi(r\delta/2) \mu\left(\Sigma_{\delta/2}(v_0)\right).
\end{align}
Splitting the integral in \eqref{posa} with $v=v_0$ into one over $\Sigma_{\delta/2}(v_0)$ and one over $\sphere\setminus\Sigma_{\delta/2}(v_0)$, and using the obvious bounds for the integrand in these cases, we get
$$
\delta\, \mu(\sphere)\le \mu\left(\Sigma_{\delta/2}(v_0)\right)+(\delta/2) \mu(\sphere)
$$
and therefore $\mu\left(\Sigma_{\delta/2}(v_0)\right)\ge (\delta/2)\mu(\sphere)$.
Substituting this into \eqref{posb}, we see that $r\le R$, where
\begin{equation}\label{sept2}
R= (2/\delta)\varphi^{-1}(2\varphi(1)/\delta),
\end{equation}
proving that $K\subset RB^n$.
\end{proof}

We can now state the first main theorem of this section, whose hypotheses allow $\psi(t)=t^p$ for $p>1$ and $G(t,u)=t^q$ for $q>0$, for example.

\bt \label{whole-general-solution-7-5}
Let $\wp: [0, \infty)\times \sphere\rightarrow [0, \infty)$ be continuous and such that $G_t>0$ on $(0, \infty)\times \sphere$ and $tG_t(t,u)$ is continuous on $[0, \infty)\times \sphere$, where $tG_t(t, u)=0$ at $t=0$ for $u\in \sphere$.  Suppose that $\psi: (0, \infty)\rightarrow (0, \infty)$ is continuous and such that $\lim_{t\rightarrow 0+} \psi(t)/t=0$ and \eqref{June-22} holds.  Then the following are equivalent:

\noindent{\rm{(i)}} The finite Borel measure $\mu$ on $\sphere$ is not concentrated on any closed hemisphere.

\noindent{\rm{(ii)}}  There exist $K\in \cK_{o}^n$ with $\inte K\neq \emptyset$ and $\tau>0$ such that $\cH^{n-1}(\Xi_K)=0$ and
\begin{equation}\label{solution-mu-7-9}
(\psi\circ h_K)\mu=\tau\,\leV(K, \cdot),
\end{equation}
where
\begin{equation}\label{july272}
\tau=\frac{1}{\leV(K, \sphere)} \int_{\sphere} \psi(h_{K}(u))\,d\mu(u).
\end{equation}
\et

\begin{proof}
Assume that (i) is true.  Following the proof of \cite[Theorem~8.2.2]{Sch}, we can construct nonzero finite discrete Borel measures $\mu_j$, $j\in \N$, such that $\mu_j\rightarrow \mu$ weakly as $j\to\infty$ and such that there is a $\delta>0$ so that (\ref{posa}) holds for $\mu$ and also with $\mu$ replaced by $\mu_j$, $j\in \N$.  In particular, $\mu_j$, $j\in \N$, is not concentrated on any closed hemisphere.  By Theorem~\ref{whole-solution-6-27}, for each $j$, there exists a convex polytope $P_j\in \cK_{(o)}^n$ such that $\mu_j=\tau_j\,\deV(P_j, \cdot)$, where
\begin{eqnarray}\label{constant-tau-7-05}
\tau_j=\frac{\mu_j(\sphere)}{\deV(P_j,\sphere)}=\frac{1}{\leV(P_j, \sphere)}
\int_{\sphere}\psi(h_{P_j}(u))\,d\mu_j(u).
\end{eqnarray}
Moreover, from (\ref{compare-ball-7-5}), we have \begin{equation}\label{compare-ball-7-5-1}
\dveV(P_j)\geq \dveV(\ball)>\int_{\sphere} G(0, u)\,du
\end{equation}
for $j\in\N$.  Theorem~\ref{whole-solution-6-27} also gives $\|h_{P_j}\|_{\varphi, \mu_j}=1$, where $\varphi$ is defined by \eqref{def-varphi-06-26}.  By Lemma~\ref{august24}, we have $P_j\subset R\ball$ for $j\in \N$, where $R$ is given by (\ref{sept2}). Then Blaschke's selection theorem implies that $P_j\rightarrow K$ for some $K\in \cK_o^n$, as $j\to \infty$, in the Hausdorff metric.  By Lemma~\ref{continuity for vG}, $\lim_{j\rightarrow \infty} \dveV(P_j)=\dveV(K)$.  This and (\ref{compare-ball-7-5-1}) imply that
\begin{equation*}
\dveV(K)\geq \dveV(\ball)>\int_{\sphere} G(0, u)\,du.
\end{equation*}
In view of (\ref{july132}), this shows that $\inte K\neq \emptyset$.

By Proposition~\ref{property-extended-7-1}(iii) and the fact that $\inte K\neq \emptyset$, we have $\leV(P_j,\cdot)\rightarrow \leV(K, \cdot)$ weakly as $j\to\infty$ and hence
$$\leV(P_j,\sphere)\to \leV(K, \sphere)>0$$
as $j\to\infty$.  Our assumption that $\lim_{t\rightarrow 0+} \psi(t)/t=0$ shows that $\psi(0)=0$ provides a continuous extension of $\psi$ to $[0, \infty)$.  This and the uniform convergence of $h_{P_j}$ to $h_K$ imply that $\psi(h_{P_j})\to \psi(h_K)$ uniformly as $j\to\infty$.  Now from the weak convergence of $\mu_j$ to $\mu$ and of $\leV(P_j,\cdot)$ to $\leV(K, \cdot)$, along with $\mu_j=\tau_j\,\deV(P_j, \cdot)$, which can be expressed in the form
$$
(\psi\circ h_{P_j})\mu_j= \tau_j\,\leV(P_j,\cdot),
$$
with $\tau_j$ as in (\ref{constant-tau-7-05}), we conclude that (\ref{solution-mu-7-9}) holds, with $\tau$ given by (\ref{july272}).  That $\tau$ is finite is a direct consequence of the continuity of $\psi$.

Since $G_t>0$, we have $\tau\geq 0$.  We claim that $\tau>0$.  To see this, use $\inte K\neq \emptyset$ to choose $v\in \sphere$ such that $\rho_K(v)>0$.
As $\mu$ is not concentrated on any closed hemisphere, the monotone convergence theorem yields
$$\lim_{j\rightarrow \infty} \int_{\Sigma_{1/j}(v)}   \langle u,  v\rangle\,d\mu(u)=\int_{\{u\in \sphere: \langle u, v\rangle> 0\}} \langle u,  v\rangle\,d\mu(u)>0,$$
where $\Sigma_{\ee}(v)$ is defined for $\ee\in (0,1)$ by (\ref{july193}). Hence a $j_0\geq 2$ exists such that
$$\mu\left(\Sigma_{1/j_0}(v)\right)\geq \int_{\Sigma_{1/j_0}(v)}\langle u,  v\rangle\,d\mu(u)>0.$$
We use this, (\ref{july272}), and the fact that
$$h_K(u)\geq \rho_K(v)\langle u, v\rangle\geq \rho_K(v)/j_0$$
for $u\in \Sigma_{1/j_0}(v)$ to obtain
\begin{eqnarray*}
\tau &\geq& \frac{1}{\leV(K, \sphere)}  \int_{\Sigma_{1/j_0}(v)} \psi(h_K(u))\,d\mu(u) \\&\geq&  \min\left\{\psi(t):t\in [\rho_K(v)/j_0 , R]\right\} \,\frac{\mu\left(\Sigma_{1/j_0}(v)\right)}{\leV(K, \sphere)}>0,
\end{eqnarray*}
proving our claim and (\ref{july272}).

It remains to be shown that $\cH^{n-1}(\Xi_K)=0$. To see this, suppose to the contrary that $\cH^{n-1}(\Xi_K)\neq 0$.  Then (see (\ref{july252})) we have $o\in \partial K$. Since $\tau>0$, we can, in view of (\ref{constant-tau-7-05}) and the fact that $\mu_j\to \mu$ and $\tau_j\to \tau$ as $j\to \infty$, assume without loss of generality that
\begin{equation}\label{july281}
\deV(P_j, \sphere)\leq \frac{2\mu(\sphere)}{\tau}<\infty
\end{equation}
for $j\in \N$, where $P_j$ is as above. Let $z\in \inte K$ be fixed.
For $E\subset \partial K$, define
$$\sigma(E)=\{z+\lambda(x-z): \ \ x\in E\ \ \mathrm{and}\ \ \lambda>0\}.
$$
Let $a, b>0$, and let $\ee>0$. From statements (a'), (b'), and (c') in the proof of \cite[Lemma~4.4]{BorFor}, we know (recall that all $P_j\in \cK_{(o)}^n$) that there exist $U\subset \partial K$ and  $j_{\ee}\in \N$ such that for $j\ge j_{\ee}$, one has $a\leq |x|\leq R$ for all $x\in \sigma(U)\cap \partial P_j$, $\cH^{n-1}(\sigma(U)\cap \partial P_j)\geq b/2$, and $h_{P_j}(u)\leq 2\ee$ if $u\in \sphere$ is an outer normal vector at $x\in \sigma(U)\cap \partial P_j$. Using these facts, (\ref{form-11-29}) for $P_j$ with $o\in\inte P_j$, and the continuity of $\wp_t$ on $(0, \infty)\times \sphere$, we obtain, for $j\geq j_{\ee}$,
\begin{eqnarray}
\deV(P_j, \sphere)&=&\frac{1}{n}\int_{\partial P_j} \frac{\langle x, \nu_{P_j}(x)\rangle}{\psi(\langle x, \nu_{P_j}(x)\rangle)}\,|x|^{1-n}\wp_t(|x|, \tilde{\pi}(x))\,dx \label{july273} \\ &\geq &  \frac{1}{n}\int_{\sigma(U)\cap \partial P_j} \frac{\langle x, \nu_{P_j}(x)\rangle}{\psi(\langle x, \nu_{P_j}(x)\rangle)}\,|x|^{1-n}\wp_t(|x|, \tilde{\pi}(x))\,dx \ge \frac{bc\,d_{\ee}}{2n},\label{uniform-bounded-dev-7-11}
\end{eqnarray}
where
$$c=\min\left\{t^{1-n}\wp_t(t, u): (t, u)\in [a, R]\times \sphere\right\}>0~~\quad{\text{and}}~~\quad d_{\ee}=\inf\{t/\psi(t): t\in (0, 2\ee]\}.$$
Since $\lim_{t\rightarrow 0+} \psi(t)/t=0$, (\ref{july281}) and (\ref{uniform-bounded-dev-7-11}) yield
$$
\infty>\frac{2\mu(\sphere)}{\tau}\geq  \deV(P_j, \sphere) \ge \lim_{\ee\to 0+} \frac{bc\, d_{\ee}}{2n}=\infty.
$$
This contradiction proves that $\cH^{n-1}(\Xi_K)= 0$.  Therefore (ii) holds.

Now assume that (ii) is true.  We claim that $\leV(K, \cdot)$ is not concentrated on any closed hemisphere; by (\ref{solution-mu-7-9}), this will yield (i).  To prove the claim, we must show that (\ref{condition for Minkowski problem}) holds when $\mu$ there is replaced by $\leV(K,\cdot)$.  If this is not true, there is a $v_0\in \sphere$ such that
\begin{equation}\label{concentration-7-8}
\int_{\sphere} \langle u,v_0\rangle_+\,d\leV(K, u)=\frac{1}{n}\int_{\sphere\cap \inte N(K,o)^*} \langle  \alpha_K(u),v_0\rangle_+ \,\rho_K(u)\,G_t(\rho_K(u), u)\,du=0,
\end{equation}
where the first equality is due to (\ref{new measue-11-27}) with $\psi\equiv 1$.  By  (\ref{july221}), we have $\rho_K(u)>0$ if $u\in \sphere \cap \inte N(K, o)^*$.  It follows from (\ref{concentration-7-8}) that $\langle \alpha_K(u), v_0\rangle_+ =0$ for ${\cH}^{n-1}$-almost all $u\in \sphere\cap \inte N(K,o)^*$. Define $X=\Xi_K\cup \sigma_K\cup Y\subset\partial K$,
where
$$Y=\left\{r_K(u)=\rho_K(u)u:  u\in \sphere\cap\inte N(K,o)^*~{\text{and}}~ \langle \alpha_K(u),v_0\rangle_+\neq 0\right\}.$$
Then the observations just made imply that $\cH^{n-1}(Y)=0$, and since $\cH^{n-1}(\Xi_K)=0$ by assumption and (\ref{zerosets}) holds, it follows that
$\cH^{n-1}(X)=0$.  Moreover, for $x=r_K(u)\in \partial K\setminus X$, we have $\langle \alpha_K(u),v_0\rangle_+=\langle \nu_K(r_K(u)),v_0\rangle_+=0$ and hence
\begin{equation}\label{august27}
\langle \nu_K(x), v_0\rangle\le 0.
\end{equation}
Next, note that if $A\subset \reg K=\partial K\setminus \sigma_K$ and $\cH^{n-1}(\partial K\setminus A)=0$, then
\begin{equation}\label{seteq}
K=\bigcap_{x\in  A}H^-(K,x),
\end{equation}
where $H^-(K,x)$ is the unique supporting halfspace of $K$ containing $K$ whose bounding hyperplane $H(K,x)$ passes through $x$.  Indeed, $K$ is contained in the set on the right-hand side of (\ref{seteq}). For the reverse inclusion, let $z\in \R^n\setminus K$. Choose a ball $B\subset\inte K$.  Then $\conv(\{z\}\cup B)\cap \partial K$ is open relative to $\partial K$ and since it has positive $\cH^{n-1}$-measure, it must contain an $x\in A$. Then $B\subset \inte H^-(K,x)$ and therefore $z\not\in H^-(K,x)$. This proves \eqref{seteq}.

The representation \eqref{seteq} immediately implies that the positive hull of $\{\nu_K(x):x\in A\}$ is $\R^n$.  Noting that $\partial K\setminus X\subset \reg K$ by the definition of $X$, we see that when $A=\partial K\setminus X$, this contradicts (\ref{august27}) and completes the proof.
\end{proof}

It is not true in general that the set $K$ in Theorem~\ref{whole-general-solution-7-5}(ii) satisfies $K\in\mathcal{K}_{(o)}^n$.  In fact this is already the situation for the $L_p$ Minkowski problem, corresponding to $G(t, u)=t^n$ and $\psi(t)=t^p$ for $p>1$; see \cite[Example~4.1]{HugLYZ}.  However, additional assumptions can be imposed ensuring that we can find a solution $K$ of \eqref{solution-mu-7-9} with $K\in\mathcal{K}_{(o)}^n$.  For example, suppose that $t/\psi(t)$ is decreasing on $(0, 1]$ and there exists $c_0>0$ such that  \begin{equation}\label{bounded-of-H-7-11}
\inf\left\{\frac{t\,\wp_t(t, u)}{\psi(t)}: (t, u)\in (0, 1]\times \sphere\right\} >nc_0.
\end{equation}
We show that it is not possible to have $o\in \partial K$ and $\inte K\neq \emptyset$.  Using (\ref{july273}), $\langle x, \nu_{P_j}(x)\rangle\leq |x|$ for $j\in \N$ and $x\in \partial P_j$, the fact that $\psi(t)/t$ is increasing on $(0, 1]$, and (\ref{bounded-of-H-7-11}), we obtain
\begin{eqnarray*}
\frac{2\mu(\sphere)}{\tau} &\geq& \frac{1}{n}\int_{\ball\cap \partial P_j} \frac{\langle x, \nu_{P_j}(x)\rangle}{\psi(\langle x, \nu_{P_j}(x)\rangle)}\,|x|^{1-n}\wp_t(|x|, \tilde{\pi}(x))\,dx\\&\geq& \frac{1}{n}\int_{\ball\cap \partial P_j} \frac{|x|}{\psi(|x|)}\,|x|^{1-n}\wp_t(|x|, \tilde{\pi}(x))\,dx \\ &\geq&  c_0\int_{\ball\cap\partial P_j} |x|^{1-n}\,dx.
\end{eqnarray*}
The argument then follows  directly from \cite[(55)-(57)]{BorFor}.  In particular, we can find $v\in \sphere$, $c_1>0$, and $0<r_0<r_1<1$ such that
\begin{equation*}
\frac{2\mu(\sphere)}{\tau} \geq  c_0\int_{\ball\cap\partial P_j} |x|^{1-n}\,dx\geq c_1\int_{B(r_1)\setminus B(r_0)} |x|^{1-n}\,dx> \frac{2\mu(\sphere)}{\tau},
\end{equation*}
where $B(r)=r\ball\cap v^{\perp}$. This contradiction proves that $K\in \cK_{(o)}^n$.

Instead of assuming the monotonicity of $\psi(t)/t$, one can assume that there exists an $\alpha\geq n-1$ such that
$$\inf\left\{t^{1-n}\,\wp_t(t, u):(t, u)\in (0, 1] \times \sphere\right\}>0 \ \ \ \mathrm{and} \ \ \
\inf_{t\in (0, 1]}\frac{t^{1+\alpha}}{\psi(t)}>0.$$
Indeed, by (\ref{july273}), we then have
\begin{equation*}
\frac{2\mu(\sphere)}{\tau}\geq  \frac{1}{n}\int_{B^n\cap \partial P_j} \frac{\langle x, \nu_{P_j}(x)\rangle}{\psi(\langle x, \nu_{P_j}(x)\rangle)}\,|x|^{1-n}\wp_t(|x|, \tilde{\pi}(x))\,dx \geq  c_2 \int_{ B^n\cap \partial P_j}   \langle x, \nu_{P_j}(x)\rangle ^{-\alpha}\,dx,
\end{equation*}
for some $c_2>0$.  It then follows directly from the arguments on \cite[p.~713]{HugLYZ} that $o\in \inte K$ and hence $K\in \cK_{(o)}^n$.

The following result provides a variant of \cite[Theorem~6.4]{GHWXY}, not requiring the condition (\ref{july195}) but with a weak additional growth condition at $0$ on $\psi$ (see the discussion after Theorem~\ref{solution-general-dual-Orlicz-main theorem-11-27}).  The hypotheses allow $\psi(t)=t^p$ for $p>0$ and $G(t,u)=t^q$ for $q<0$, for example.

\bt\label{whole-general-solution-7-9-decreaing}
Let $\wp:(0, \infty)\times \sphere\rightarrow (0, \infty)$ be continuous and such that $\wp_t$ is continuous and negative on $(0, \infty)\times \sphere$.  Let $0<\ee_0<1$ and suppose that \eqref{condE2-6-30} holds for $v\in \sphere$.
Suppose that $\psi: (0, \infty)\rightarrow (0, \infty)$ is continuous, \eqref{June-22} holds, and that $\varphi$ is finite when defined by \eqref{def-varphi-06-26}.  Then the following are equivalent:

\noindent{\rm{(i)}} The  finite Borel measure $\mu$ on $\sphere$ is not concentrated on any closed hemisphere.

\noindent{\rm{(ii)}}  There exist $K\in \cK_{(o)}^n$  and $\tau<0$  such that  \begin{equation}\label{july241}
\mu=\tau\,\deV(K,\cdot),
\end{equation}
where
\begin{equation}\label{july242}
\tau=\frac{\mu(\sphere)}{\deV(K, \sphere)}.
\end{equation}
\et

\begin{proof}  Assume that (i) is true. Define $\varphi$ as in \eqref{def-varphi-06-26}.  As at the beginning of the proof of Theorem~\ref{whole-general-solution-7-5}, but using Theorem~\ref{solution-general-dual-Orlicz-main theorem-11-27} instead of Theorem~\ref{whole-solution-6-27}, we can find nonzero finite discrete Borel measures $\mu_j$, $j\in \N$, not concentrated on any closed hemisphere, such that $\mu_j\rightarrow \mu$ weakly as $j\to\infty$, and convex polytopes $P_j\in \cK_{(o)}^n$ such that $\mu_j=\tau_j \,\deV(P_j, \cdot)$, where (in view of (\ref{solution-explicit-6-30})) $\tau_j$ satisfies (\ref{constant-tau-7-05}) and $\|h_{P_j}\|_{\varphi,\mu_j}=1$ for $j\in \N$.  From the latter property and Lemma~\ref{august24}, it follows as in the proof of Theorem~\ref{whole-general-solution-7-5} that $(P_j)_{j\in\N}$ is bounded. Hence, we can extract a subsequence that converges to $K\in\K^n_o$. Next, we show that $o\in \inte K$. In fact, if  $o\in \partial K$, we can apply Lemma~\ref{decreasing-general-6-30} to get $\lim_{j\rightarrow \infty} \dveV(P_j)=\infty$.  However, since $P_j$ corresponds to $P(z^0)$ in Theorem~\ref{solution-general-dual-Orlicz-main theorem-11-27}, (\ref{compare-ball-7-9}) implies that
\begin{eqnarray*}
\dveV(P_j) \le  \dveV(\ball)<\infty
\end{eqnarray*}
for all $j\in \N$, a contradiction.

Then (\ref{july241}) and (\ref{july242}) follow from the weak convergence of $\mu_j$ to $\mu$ and of $\deV(P_j, \cdot)$ to $\deV(K, \cdot)$, the latter a consequence of \cite[Proposition~6.2(ii)]{GHWXY}. In particular, we use that $\deV(P_j, \sphere)\to \deV(K, \sphere)\in (0,\infty)$ to ensure the convergence of $(\tau_j)_{j\in\N}$.

Suppose that (ii) holds.  By \cite[Proposition~6.2(iii)]{GHWXY}, $\deV(K, \cdot)$ is not concentrated on any closed hemisphere, so by (\ref{july241}), this is also the case for $\mu$.
\end{proof}

The final result in this section addresses the uniqueness problem related to Theorem~\ref{whole-general-solution-7-9-decreaing} and generalizes and extends \cite[Theorem~8.3]{LYZ-Lp}.  It can be applied, for example, when $G(t,u)=t^q$, $q\neq 0$, and $\psi(s)=s^p$ with $q<p$. Note that when $\psi\equiv 1$ and $G_t<0$, the result holds for general $K,K'\in \cK^n_{(o)}$ by Theorem~\ref{theorem5.4new}, since the assumption there that $tG_t(t,u)$ is strictly increasing in $t$ for $u\in\sphere$ implies the second inequality in (\ref{condna}).  We do not know if the result holds for general $\psi$ and general $K,K'\in \cK^n_{(o)}$.

\bt\label{uniquemp}
Let $\wp: (0, \infty)\times \sphere\rightarrow (0, \infty)$ and $\psi:(0,\infty)\rightarrow (0, \infty)$ be continuous.  Suppose that $G_t>0$ (or $G_t<0$) on $(0,\infty)\times\sphere$ and that if
\begin{equation}\label{condna}
\frac{ G_t( t,u)}{\psi( s)}\ge \frac{\lambda G_t(\lambda t,u)}{\psi(\lambda s)}~\quad{\text{(or}}~\quad \frac{ G_t( t,u)}{\psi( s)}\le \frac{\lambda G_t(\lambda t,u)}{\psi(\lambda s)},~~{\text{respectively)}}
\end{equation}
for some $\lambda, s,t>0$ and $u\in\sphere$, then $\lambda\ge 1$. If $K,K'\in \cK^n_{(o)}$ are both polytopes or both have support functions in $C^2$ and $\deV(K, \cdot)=\deV(K',\cdot)$, then $K=K'$.
\et

\begin{proof}
Suppose that $K,K'\in \cK^n_{(o)}$ are such that $\deV(K, \cdot)=\deV(K',\cdot)$ and $K\neq K'$.  Then we can assume without loss of generality that $K\not\subset K'$, so there is a maximal $\lambda<1$ such that $\lambda K\subset K'$.

Consider first the case when $K$ and $K'$ are polytopes.  By Lemma~\ref{facet-06-27}, the facets of $K$ and $K'$ have the same outer unit normal vectors, $u_1,\dots,u_m$, say, and from (\ref{polytope-06-25}) and (\ref{july152}), we have
$$\deV(K, \cdot)=\deV(K',\cdot)=\sum_{i=1}^m\gamma_i\delta_{u_i},$$
where
\begin{equation}\label{sept251}
\gamma_i= \int_{\tilde{\pi}(F(K, u_i))}\frac{\rho_{K}(u)\, \wp_t(\rho_K(u),u)}{n\psi(h_{K}(u_i))} \,du= \int_{\tilde{\pi}(F(K', u_i))} \frac{\rho_{K'}(u)\, \wp_t(\rho_{K'}(u),u)}{n\psi(h_{K'}(u_i))}\,du.
\end{equation}
Since the facets of $\lambda K$ and $K'$ also have the same outer unit normal vectors and $\lambda$ is maximal, at least one facet of $\lambda K$ is contained in a facet of $K'$.  If this facet has outer unit normal vector $u_i$, then
\begin{equation}\label{sept252}
h_{\lambda K}(u_i)=h_{K'}(u_i),~~\quad~~\tilde{\pi}(F(K, u_i))=\tilde{\pi}(F(\lambda K, u_i))\subset\tilde{\pi}(F(K', u_i)),
\end{equation}
and
\begin{equation}\label{sept253}
\rho_{\lambda K}(u)=\rho_{K'}(u)~~\quad~~{\text{for}}~~\quad~~ u\in \tilde{\pi}(F(K, u_i)).
\end{equation}
If $G_t>0$ (the argument when $G_t<0$ is similar), we conclude from (\ref{sept251}), (\ref{sept252}), and (\ref{sept253}) that
$$
\int_{\tilde{\pi}(F(K, u_i))}\frac{\rho_{K}(u)\, \wp_t(\rho_K(u),u)}{n\psi(h_{K}(u_i))}\,du
\ge  \int_{\tilde{\pi}(F(K, u_i))}\frac{\rho_{\lambda K}(u)\, \wp_t(\rho_{\lambda K}(u),u)}{n\psi(h_{\lambda K}(u_i))} \,du.
$$
Since $\mathcal{H}^{n-1}(\tilde{\pi}(F(K, u_i))>0$, there is a $u\in \tilde{\pi}(F(K, u_i))$ such that
$$\frac{\rho_{K}(u)\, \wp_t(\rho_K(u),u)}{n\psi(h_{K}(u_i))}\ge
\frac{\rho_{\lambda K}(u)\, \wp_t(\rho_{\lambda K}(u),u)}{n\psi( h_{\lambda K}(u_i))},$$
that is,
$$\frac{\wp_t(\rho_K(u),u)}{\psi(h_{K}(u_i))}\ge
\frac{\lambda \, \wp_t(\lambda \rho_{K}(u),u)}{\psi(\lambda h_{K}(u_i))}.$$
Now the first inequality in (\ref{condna}) with $s=h_K(u_i)$ and $t=\rho_K(u)$ yields $\lambda\ge 1$, a contradiction.  This completes the proof for when $K$ and $K'$ are polytopes.

For the other case, note firstly that if $L\in\cK^n_{(o)}$ and $h_L\in C^2$, then $S(L,\cdot)$ is absolutely continuous with respect to $\mathcal{H}^{n-1}$
with continuous density $R(L,\cdot)$, where $R(L,u)$ is the product of the principal radii of curvature of $L$ at $u\in\sphere$.  (This is well known when $L$ is of class $C^2_+$; see, for example, \cite[(4.26), p.~217]{Sch}.  When $h_L\in C^2$, one can observe that \cite[Lemma~5.1]{BFH} implies that \cite[Theorem~3.7(c)]{Hug2} holds, and then \cite[Theorem~3.7(a)]{Hug2} yields the absolute continuity of $S(L,\cdot)$.  The form of the density is then given by \cite[Theorem~3.5]{Hug1}.)  Let $K,K'\in \cK^n_{(o)}$ and $h_K, h_{K'}\in C^2$. Using \cite[(24)]{GHWXY}, we obtain
\begin{eqnarray}
&&\frac{h_K(u)}{\psi (h_K(u))}\,|\nabla h_K(u)|^{1-n}\,
\wp_t\left(|\nabla h_K(u)|, \nabla h_K(u)/|\nabla h_K(u)|\right)R(K,u)\nonumber\\
&&\qquad =\frac{h_{K'}(u)}
{\psi (h_{K'}(u))}\,|\nabla h_{K'}(u)|^{1-n}\,
\wp_t\left(|\nabla h_{K'}(u)|, \nabla h_{K'}(u)/|\nabla h_{K'}(u)|\right)R(K',u)\label{condnb}
\end{eqnarray}
for all $u\in\sphere$, since both sides of (\ref{condnb}) are continuous functions.  Since $\lambda K\subset K'$ and $\lambda<1$ is maximal, there exists $u_0\in\sphere$ such that $h_{\lambda K}(u_0)=h_{K'}(u_0)$ and $\nabla h_{\lambda K}(u_0)=\nabla h_{K'}(u_0)$, i.e., $\lambda K$ and $K'$ have a common boundary point with common outer unit normal vector $u_0$.

We claim that
\begin{equation}\label{sept23}
R(K',u_0)\ge R(\lambda K,u_0)=\lambda^{n-1}R(K,u_0).
\end{equation}
It suffices to prove the inequality, since the equality follows by homogeneity.  Let $u=u_0+av$, where $a>0$ and $v\in\sphere$.  For $L\in\cK^n_{(o)}$ with $h_L\in C^2$, and $u\in \sphere$, let $d^2 h_L[u]$ denote the second differential of $h_L$ at $u$, considered as a bilinear form on $\R^n$.  Since $h_{\lambda K}\le h_{K'}$, $h_{\lambda K}(u_0)=h_{K'}(u_0)$, and $\nabla h_{\lambda K}(u_0)=\nabla h_{K'}(u_0)$, we may apply the first displayed equation in \cite[p.~31, Note~3]{Sch} (with $f=h_{L}$, $Af(x)=d^2 h_L[x]$, $x=u_0$, and $y=u$, for $L=\lambda K$ and $L=K'$), to obtain
$$
\frac{1}{2}d^2 h_{\lambda K}[u_0](av,av)+r_{\lambda K}(u_0,a)a^2\le \frac{1}{2}d^2 h_{K'}[u_0](av,av)+r_{K'}(u_0,a)a^2,
$$
where $r_{\lambda K}(u_0,a),r_{K'}(u_0,a)\to 0$ as $a\to 0+$. Dividing by $a^2$ letting $a\to 0+$, we get
$d^2 h_{\lambda K}[u_0](v, v)\le d^2h_{K'}[u_0](v, v)$ for $v\in \sphere$. We write $d^2 h_{\lambda K}[u_0]|u_0^\perp$ and $d^2 h_{K'}[u_0]|u_0^\perp$ for the symmetric, positive semidefinite linear maps from $u_0^\perp$ to itself, associated with the restrictions of the bilinear forms to $u_0^\perp\times u_0^\perp$. By \cite[Corollary~2.5.2]{Sch}, which in particular guarantees that for both maps $u_0$ is an eigenvector with eigenvalue zero, \cite[p.~124, l.~-3]{Sch}, and with the help of \cite[Corollary~7.7.4(e)]{HornJohnson}, we conclude that
$$
R(\lambda K,u_0)=\det\left(d^2 h_{\lambda K}[u_0]|u_0^\perp\right)\le \det\left(d^2 h_{K'}[u_0]|u_0^\perp\right)=R(K',u_0),
$$
proving the claim.

Suppose that $G_t>0$ on $(0,\infty)\times\sphere$; a similar argument applies when $G_t<0$ instead. By \eqref{condnb} with $u=u_0$, and (\ref{sept23}), we have
\begin{eqnarray*}
\lefteqn{\frac{h_K(u_0)}{\psi (h_K(u_0))}\,|\nabla h_K(u_0)|^{1-n}\,
\wp_t\left(|\nabla h_K(u_0)|, \nabla h_K(u_0)/|\nabla h_K(u_0)|\right)R(K,u_0)}\\
&=&\frac{h_{\lambda K}(u_0)}{\psi (h_{\lambda K}(u_0))}\,|\nabla h_{\lambda K}(u_0)|^{1-n}\,
\wp_t\left(|\nabla h_{\lambda K}(u_0)|, \nabla h_{\lambda K}(u_0)/|\nabla h_{\lambda K}(u_0)|\right)R(K',u_0)\\
&\ge&\lambda \frac{h_K(u_0)}{\psi (\lambda h_{  K}(u_0))}\,\lambda^{1-n} |\nabla h_{ K}(u_0)|^{1-n}\,
\wp_t\left(\lambda|\nabla h_{ K}(u_0)|, \nabla h_{K}(u_0)|/\nabla h_{K}(u_0)|\right)\lambda^{n-1}R(K,u_0).
\end{eqnarray*}
Therefore
$$
\frac{\wp_t\left(|\nabla h_K(u_0)|, \nabla h_K(u_0)/|\nabla h_K(u_0)|\right)}{\psi (h_K(u_0))}\ge
\frac{\lambda \wp_t\left(\lambda |\nabla h_{ K}(u_0)|, \nabla h_{  K}(u_0)/|\nabla h_{  K}(u_0)|\right)}{\psi (\lambda h_K(u_0))}.
$$
But then the first inequality in \eqref{condna} implies that $\lambda\ge 1$, a contradiction proving that $K=K'$.
\end{proof}

\section{Minkowski problems for even measures}\label{section-7}
In this section we revisit the Minkowski problems considered in earlier sections, focusing on the case of even measures and attempting to keep the discussion as brief as possible. Recall that $\cK_{os}^n$ (or $\cK_{(o)s}^n$) denote the class of origin-symmetric compact convex sets containing the origin (or containing the origin in their interiors, respectively).  Also, note that an even measure is not concentrated on any closed hemisphere if and only if it is not concentrated on a great subsphere.

The hypotheses of the next theorem allow $\psi(t)=t^p$ for $p>0$ and $G(t,u)=t^q$ for $q>0$, for example.

\bt \label{whole-general-solution-7-55-even}
Let $\wp: [0, \infty)\times \sphere\rightarrow [0, \infty)$ be continuous and such that $G_t$ is continuous and positive on $(0, \infty)\times \sphere$.
Assume that $G_t(t,u)=G_t(t,-u)$ for $(t,u)\in (0,\infty)\times\sphere$. Suppose that $\psi: (0, \infty)\rightarrow (0, \infty)$ is continuous, \eqref{June-22} holds, and that $\varphi$ is finite when defined by \eqref{def-varphi-06-26}.    Then the following are equivalent:

\noindent{\rm{(i)}} The  finite even Borel measure $\mu$ on $\sphere$ is not concentrated on any closed hemisphere.

\noindent{\rm{(ii)}}  There is a $K\in \cK_{(o)s}^n$  such that $\mu=\tau \deV(K,\cdot)$, with $\tau>0$ as in  \eqref{july272}.
\et

\begin{proof}
We first observe that under our extra assumption that $G_t(t,u)=G_t(t,-u)$ for $(t,u)\in (0,\infty)\times\sphere$, Theorem~\ref{whole-solution-6-27} holds for even discrete measures.  Specifically, if
$$
\mu=\sum_{i=1}^m \lambda_i(\delta_{u_i}+\delta_{-u_i}),
$$
where $\lambda_i>0$ for $i=1,\dots,m$ and $\{\pm u_1, \pm u_2,\dots,\pm u_m \}\subset \sphere$, there is a convex polytope $P\in {\cK}^n_{(o)s}$ satisfying (\ref{july282}).  Indeed, the proof of Theorem~\ref{whole-solution-6-27} can be easily adapted, as follows.  For each $z=(z_1, \dots, z_m)\in [0,\infty)^m$, let
$$
P_e(z) =\{x\in\R^n:\,|\langle x, u_i\rangle|\leq z_i, \ \mbox{for}\ i=1,\dots, m\},
$$
so that $P_e(z)$ is a convex polytope in $\cK_{os}^n$.  As in the proof of  Theorem~\ref{whole-solution-6-27}, one can find $z^0=(z_1^0, \dots, z_m^0)\in M_+$ such that
$$
\dveV(P_e(z^0))= \max \left\{\dveV(P_e(z)): z\in M_+\right\}.
$$
Moreover, (\ref{compare-ball-7-5}) holds with $P(z^0)$ replaced by $P_e(z^0)$.  From this, we see that $P_e(z^0)\in\cK_{(o)s}^n$ and $z_i^0>0$ for $i=1, \dots, m$.  One can adjust the argument used to prove (\ref{variation-11-27-1-1}) to obtain
\begin{equation}\label{symmetric formula}
\frac{\partial \dveV(P_e(z))}{\partial z_i}\bigg|_{z=z^0}= 2n \frac{\leV(P_e(z^0), \{u_i\})}{h_{P_e(z^0)}(u_i)}
\end{equation}
for $i=1,\dots,m$. The method of Lagrange multipliers provides $\tau\in \R$ such that
\begin{equation}\label{july284}
\frac{ \tau}{2n}\frac{\partial \dveV(P_e(z))}{\partial z_i}\bigg|_{z=z^0} =\frac{\partial \sum_{i=1}^m\lambda_i \varphi(z_i)}{\partial z_i}\bigg|_{z=z^0}
\end{equation}
for $i=1, \dots, m$.  Then (\ref{symmetric formula}) and (\ref{july284}) can be used instead of (\ref{variation-11-27-1-1}) and (\ref{july171}), respectively, and the rest of the proof of Theorem~\ref{whole-solution-6-27} can be followed up to \eqref{solution-explicit-6-30} to conclude the proof in the case of an even discrete measure.

With Theorem~\ref{whole-solution-6-27} for even discrete measures in hand, the proof of (i)$\Rightarrow$(ii) in Theorem~\ref{whole-general-solution-7-5} can be followed without difficulty to obtain the same implication for even measures, where $K$ is origin symmetric. In particular, we can take advantage of the fact that it easily follows that $P_j\to K\in \K^n_{os}$ and $\inte K\neq \emptyset$, hence $K\in\K^n_{(o)s}$. But then $h_K$ is bounded away from zero and no continuous extension of $\psi$ at $0$ is needed.

The implication (ii)$\Rightarrow$(i) follows from the proof of the same implication in Theorem~\ref{whole-general-solution-7-5} together with the evenness of $\leV(K,\cdot)$ when $K$ is origin symmetric and our extra assumption on $G$ holds. (Recall that $\Xi_K=\emptyset$ if $K\in\K^n_{(o)}$.)
\end{proof}

We omit the proof of the following result, which provides the even analogue of \cite[Theorem~6.4]{GHWXY}, since it follows without difficulty from the argument given in the proof of Theorem~\ref{whole-general-solution-7-9-decreaing}.  The hypotheses allow $\psi(t)=t^p$ for $p>0$ and $G(t,u)=t^q$ for $q<0$, for example.

\bt \label{whole-general-solution-7-9-decreaing-even}
Let $\wp: (0, \infty)\times \sphere\rightarrow (0, \infty)$ and $\psi: (0, \infty)\rightarrow (0, \infty)$ satisfy the assumptions of Theorem~\ref{whole-general-solution-7-9-decreaing} and suppose also that $G_t(t,u)=G_t(t,-u)$ for $(t,u)\in (0,\infty)\times\sphere$.  Then Theorem~\ref{whole-general-solution-7-9-decreaing} holds when in {\rm{(i)}} $\mu$ is an even measure and in {\rm{(ii)}} $K$ is origin symmetric.
\et

Our final result addresses Problem~\ref{Minkowski-c-06-25} when $G_t<0$ and $\psi$ is decreasing.  The hypotheses allow $\psi(t)=t^p$ for $p<0$ and $G(t,u)=t^q$ for $q<0$, for example.

If $\psi: (0, \infty)\rightarrow (0, \infty)$ is continuous, define
\begin{equation}\label{decreasing-varphi-7-27} \overline{\varphi}(t)=\int_t^{\infty}\frac{\psi(s)}{s}\,ds
\end{equation}
for $t>0$.

\bt \label{solution-general-dual-Orlicz-main theorem-11-27-even}
Let $\mu$ be a nonzero finite even Borel measure vanishing on great subspheres.  Let $\wp$ and $\wp_t$ be continuous on $(0, \infty)\times \sphere$, where $\wp_t<0$ and where $\wp_t(t, u)=\wp_t(t, -u)$ for $(t,u)\in (0,\infty)\times\sphere$. Suppose that there is some  $0<\ee_0<1$ such that \eqref{condE2-6-30} holds for $v\in \sphere$.  Let $\psi: (0, \infty)\rightarrow (0, \infty)$ be continuous and suppose that $\overline{\varphi}$ is finite when defined by \eqref{decreasing-varphi-7-27}. Then there exists a $K \in \cK_{(o)s}^n$ such that
\begin{equation}\label{msol-even}
\frac{\mu}{\mu(\sphere)}=\frac{\deV(K, \cdot)}{\deV(K, S^{n-1})}.
\end{equation}
\et

\begin{proof}
Since $\wp_t<0$, we may define $a_0\in [0,\infty)$ by
$$
a_0=\lim_{t\to \infty} \int_{\sphere} \wp(t,u)\,du.
$$
Define the functional $F: C^+(\sphere)\rightarrow \R$ by
$$
F(f)=\frac{1}{\mu(\sphere)+a_0}\int_{S^{n-1}}\overline{\varphi}
(f(u))\,d\mu(u)
$$
for $f\in C^+(\sphere)$, and define $F(K)=F(h_K)$ for $K \in \cK_{(o)s}^n$.  Let
\begin{equation}\label{optimization-convex-body-11-28-11-even}
\alpha=\sup\left\{F(K): K\in \cK_{(o)s}\text{ and } \dveV(K)=\mu(\sphere)+a_0\right\}.
\end{equation}
As in the proof of \cite[Theorem~6.4]{GHWXY}, there is an $r_0>0$ such that $\dveV(r_0B^n)=\mu(\sphere)+a_0$, so the supremum in (\ref{optimization-convex-body-11-28-11-even}) is taken over a nonempty set. (Note that our assumptions on the even measure $\mu$ imply in particular that it is not concentrated on any closed hemisphere, as is assumed in \cite[Theorem~6.4]{GHWXY}.)  Choose $K_j\in\cK_{(o)s}^n$, $j\in\N$, such that $\dveV(K_j)=\mu(\sphere)+a_0$ and $\lim_{j\rightarrow \infty} F(K_j)=\alpha$. The proof of \cite[Theorem~6.4]{GHWXY} shows that there is an $R>0$ such that the polar bodies satisfy $K_j^{*} \subset R \ball$ for $j\in \N$.  By relabeling, if necessary, using Blaschke's selection theorem, and noting that $K_j^*$ is also origin symmetric for $j\in \N$, a $Q_0\in \cK_{os}^n$ can be found such that $K_j^*\rightarrow Q_0$ as $j\to\infty$.

Define $\widetilde{\varphi}$ by $\widetilde{\varphi}(t)=\overline{\varphi}(1/t)$ for $t>0$. The dominated convergence theorem shows
that $\overline{\varphi}(t)\to 0$ as $t\to\infty$. Then our assumption on $\overline{\varphi}$ implies that
$$
\widetilde{\varphi}(0)=\lim_{t\rightarrow 0+} \widetilde{\varphi}(t)=\lim_{t\rightarrow 0+}  \overline{\varphi} (1/t)=0
$$
defines a continuous extension of $\widetilde{\varphi}$ at $0$.
By (\ref{bi-polar--1}), we have
$$
F(h_{K_j})=\frac{1}{\mu(\sphere)+a_0}\int_{S^{n-1}}\overline{\varphi}(h_{K_j}(u))
\,d\mu(u)
=\frac{1}{\mu(\sphere)+a_0}\int_{S^{n-1}}\widetilde{\varphi}(\rho_{K^*_j}(u))\,d\mu(u).
$$
We claim that $Q_0\in \cK_{(o)s}^n$.  In fact, assume that $\inte Q_0=\emptyset$, so that $Q_0\subset v ^{\perp}$ for some $v \subset\sphere$. Then, as shown in the proof of Lemma~\ref{continuity for vG}, $\rho_{K_j^*}(u)\to 0$ as $j\to\infty$ for $u\in \sphere \setminus v^\perp$. Since $\widetilde{\varphi}:[0,\infty)\to[0,\infty)$ is continuous, it follows that $\widetilde{\varphi}(\rho_{K_j^*}(u))\to \widetilde{\varphi}(0)=0$ as $j\to\infty$ for  $u\in\sphere\setminus v^\perp$, and hence for $\mu$-almost all $u\in\sphere$, as $\mu$ vanishes on the great subsphere $\sphere\cap v^\perp$.  The continuity of $\widetilde{\varphi}$ also implies that $M_1=\max\{\widetilde{\varphi}(t): t\in [0, R]\}<\infty$.  Hence the dominated convergence theorem can be applied and yields
$$
\alpha=\lim_{j\to\infty} F(h_{K_j}) =  \lim_{j\to\infty} \frac{1}{\mu(\sphere)+a_0}\int_{S^{n-1}}\widetilde{\varphi}
(\rho_{K^*_j}(u))\,d\mu(u) =0.
$$
But this is impossible because $\alpha\geq F(r_0\ball)=\overline{\varphi}(r_0)>0$. This proves the claim.

Let $K_0=Q_0^*$.  Then $K_0\in\cK_{(o)s}^n$.  Also, $K_j\rightarrow K_0$ as $j\to\infty$, so  $\dveV(K_j)\rightarrow \dveV(K_0)$ as $j\to\infty$ by the continuity of $\wp$ and \cite[Lemma~6.1]{GHWXY}, yielding $\dveV(K_0)=\mu(\sphere)+a_0$.   If $f\in C^+(\sphere)$, the support function $h_{[f]}$ of the Wulff shape $[f]$ of $f$, defined by (\ref{july287}), satisfies $h_{[f]}\leq f$.  As $\overline{\varphi}$ is decreasing, we have $\overline{\varphi}(h_{[f]})\geq \overline{\varphi}(f)$. Consequently,
\begin{equation}\label{sept83}
F(h_{K_0})=\alpha=\sup\left\{F(f): \dveV([f])=\mu(\sphere)+a_0 \ \mathrm{and}\ f\in C^+(\sphere) \ \ \mathrm{is\ even}\right\}.
\end{equation}
Let $g\in C(\sphere)$ be even.  We apply the method of Lagrange multipliers, following the argument in the proof of \cite[Theorem~6.4]{GHWXY} from \cite[(71)]{GHWXY} onwards, where (\ref{variation-11-27-1}) and $\overline{\varphi}$ play the role of \cite[(59)]{GHWXY} and $\varphi$.  (Note that from (\ref{decreasing-varphi-7-27}), we have
$\psi(t)=-t\overline{\varphi}'(t)$.)  The extra constant $a_0$ in (\ref{sept83}) has no effect on the conclusion, which is that
$$
\int_{S^{n-1}} g(u)\,d\mu(u) = -n\tau\int_{\sphere} g(u)\,d\deV(K_0, u),
$$
where
$$
\tau = -\frac{\mu(\sphere)}{n\,\deV(K_0,\sphere)}.
$$
As $g$ is an arbitrary even function in $C(\sphere)$, we can use our assumption that $G_t(t,u)=G_t(t,-u)$ for $(t,u)\in (0,\infty)\times\sphere$ to obtain (\ref{msol-even}) with $K$ replaced by $K_0$.
\end{proof}

	\end{document}